\documentclass[11pt, reqno]{amsart}
\usepackage{mathrsfs}
\usepackage[active]{srcltx}
\usepackage{mathrsfs,amsmath}
\usepackage{mathtools}
\usepackage{longtable}
\usepackage{todonotes}
\usepackage[autostyle]{csquotes}
\allowdisplaybreaks

\UseRawInputEncoding

\usepackage{xcolor}
\definecolor{bluecite}{HTML}{0875b7}
\usepackage[unicode=true,
bookmarksopen={true},
pdffitwindow=true,
colorlinks=true,
linkcolor=bluecite,
citecolor=bluecite,
urlcolor=bluecite,
hyperfootnotes=false,
pdfstartview={FitH},
pdfpagemode= UseNone]{hyperref}

\newtheorem{proposition}{Proposition}[section]
\newtheorem{theorem}{Theorem}[section]

\newtheorem{remark}{Remark}[section]
\numberwithin{equation}{section}

\usepackage{bbm}  
\usepackage{dsfont}

\address{\textsc{Alexandru Krist\'aly}: Department of Economics, Babe\c s-Bolyai University, str. Teodor Mihali 58-60, 400591, Cluj-Napoca, Romania \& 
	Institute of Applied Mathematics, \'Obuda
	University, 
	Budapest, Hungary.} 
\email{alex.kristaly@econ.ubbcluj.ro; kristaly.alexandru@uni-obuda.hu}

\thanks{Research supported by the
	Excellence Researcher Program \'OE-KP-2-2022 of \'Obuda University, Hungary.
}

\title[Sharp Sobolev inequalities on Riemannian manifolds]{Sharp Sobolev inequalities on noncompact Riemannian manifolds with ${\sf Ric}\geq 0$ via Optimal Transport theory}

\author[Alexandru Krist\'aly]{Alexandru Krist\'aly}

\begin{document}
%	\vspace{-1cm}
	\begin{abstract}
%		In their seminal work, Cordero-Erausquin,  Nazaret and Villani [\textit{Adv.\ Math.}, 2004] established the sharp Sobolev inequalities in Euclidean spaces (endowed with arbitrary norm) via \textit{optimal mass transport theory}, raising the question whether their approach is  powerful enough to yield sharp Sobolev inequalities also on Riemannian manifolds. In the present paper we affirmatively answer their question; namely, by using exclusively optimal mass transportation, sharp $L^p$-Sobolev and $L^p$-logarithmic Sobolev inequalities (both for $p>1$ and $p=1$) are established on Riemannian manifolds of non-negative Ricci curvature,  where the optimal constants contain the \textit{asymptotic volume growth}. As a byproduct of our approach, we give an alternative, elementary proof to the main result of do Carmo and Xia [\textit{Comp.\ Math.}, 2004] concerning the quantitative volume non-collapsing estimates under the validity of certain Sobolev inequalities. 
		In their seminal work, Cordero-Erausquin,  Nazaret and Villani [\textit{Adv.\ Math.}, 2004] proved sharp Sobolev inequalities in Euclidean spaces via \textit{Optimal  Transport}, raising the question whether their approach is  powerful enough to produce sharp Sobolev inequalities also on Riemannian manifolds. By using $L^1$-optimal transport  approach, the compact case  has been successfully treated by Cavalletti and Mondino [\textit{Geom.\ Topol.}, 2017], even on metric measure spaces verifying the synthetic lower Ricci	curvature bound. 
		In the present paper we affirmatively answer the above question for noncompact Riemannian manifolds with non-negative Ricci curvature; namely, by using Optimal  Transport theory with quadratic distance cost, sharp $L^p$-Sobolev and $L^p$-logarithmic Sobolev inequalities (both for $p>1$ and $p=1$) are established,  where the sharp constants contain the \textit{asymptotic volume ratio} arising from precise asymptotic properties of the Talentian and Gaussian bubbles, respectively. As a byproduct, we give an alternative, elementary proof to the main result of do Carmo and Xia [\textit{Compos.\ Math.}, 2004] and subsequent results,	concerning the quantitative volume non-collapsing estimates on Riemannian manifolds with non-negative Ricci curvature that support Sobolev inequalities.

%				Some rigidity results are also derived in the context of Riemannian manifolds with nonnegative Ricci curvature.  
	\end{abstract}
	\maketitle 

%	\tableofcontents
\begin{center}
\textit{Dedicated to Professor L\'aszl\'o Lov\'asz on the occasion of his 75th anniversary}
\end{center}
	%	\vspace{0.7cm}
	\section{Introduction and Main Results}
	
The primary goal of the present paper is to answer a question of Cordero-Erausquin,  Nazaret and Villani \cite{CE-N-Villani}; we show that the \textit{Optimal  Transport} (OT, for short) for quadratic distance cost functions can be efficiently applied  to establish sharp Sobolev inequalities on noncompact, complete   Riemannian manifolds with non-negative Ricci curvature (${\sf Ric}\geq 0$, for short).
%, the optimal Sobolev constants containing the \textit{asymptotic volume ratio}. 
%We perform our arguments for the $L^p$-Sobolev and $L^p$-logarithmic Sobolev inequalities (both for $p>1$ and $p=1$). 

After Aubin \cite{Aubin} initiated the famous  AB-program,  considerable efforts have been made to determine the optimal constants in Sobolev inequalities, both in euclidean and non-euclidean settings; see also Druet and Hebey \cite{Druet-Hebey}  and Hebey \cite{Hebey}. 
A central place within this program is occupied by the 
%In conformity to the question in \cite{CE-N-Villani}, we are going to use OT theory for quadratic cost functions, developed by McCann \cite{McCann}, combined with fine asymptotic properties of the Talentian/Guassian bubbles on Riemannian manifolds with non-positive Ricci curvature. 
classical  $L^p$-Sobolev inequality in $\mathbb R^n$;  given $n\geq 2$, $p\in (1,n)$ and the critical Sobolev exponent $p^\star=\frac{pn}{n-p}$, the latter inequality states that
	\begin{equation}\label{Sobolev-000}
	\displaystyle  { \left(\int_{\mathbb R^n} |f|^{p^\star}{\rm d}x\right)^{1/p^\star}\leq
		{\sf AT}(n,p)
		\left(\int_{\mathbb R^n} |\nabla f|^p {\rm d}x\right)^{1/p}},\ \ \forall  f\in C_0^\infty(\mathbb R^n),
\end{equation}
where 
$${\sf AT}(n,p)=\pi^{-\frac{1}{2}}n^{-\frac{1}{p}} \left(\frac{p-1}{n-p}\right)^{1/p'}\left(\frac{\Gamma(1+n/2)\Gamma(n)}{\Gamma(n/p)\Gamma(1+n/p')}
\right)^{{1}/{n}},$$ and $p'=\frac{p}{p-1}$ is the conjugate of $p$, see Aubin \cite{Aubin} and Talenti \cite{Talenti}.
In addition, the Aubin--Talenti-constant 	${\sf AT}(n,p)$ in  \eqref{Sobolev-000} is sharp,  which is achieved by the class $(f_\lambda)_{\lambda>0}$  of  Talentian bubbles
\begin{equation}\label{talentian-bubble-eukl}
	f_\lambda(x)=(\lambda+|x|^{p'})^\frac{p-n}{p},\ x\in \mathbb R^n.
\end{equation}
In the limit case $p=1$, we have
${\sf AT}(n,1)=\lim_{p\searrow 1}{\sf AT}(n,p)=n^{-1}\omega_n^{-\frac{1}{n}},$ where $\omega_n=\frac{\pi^{n/2}}{\Gamma(\frac{n}{2}+1)}$ is the volume of the unit ball in $\mathbb R^n$, while \eqref{Sobolev-000} reduces to the sharp $L^1$-Sobolev inequality, which is in turn equivalent to the sharp isoperimetric inequality in $\mathbb R^n$. To prove \eqref{Sobolev-000},  Aubin \cite{Aubin} and Talenti \cite{Talenti} stated the P\'olya-Szeg\H o inequality by using  Schwarz-symmetrization and the sharp isoperimetric inequality in $\mathbb R^n$, reducing a multidimensional problem to a more tractable $1$-dimensional one. 

A genuinely new proof for the sharp Sobolev inequality \eqref{Sobolev-000} has been obtained by Cordero-Erausquin,  Nazaret and Villani \cite{CE-N-Villani} by using  $L^2$-OT-theory (i.e., for the quadratic distance cost function), based on the Brenier map \cite{Brenier}. In their approach the linear structure of the Euclidean space $\mathbb R^n$ played a crucial role. Regarding the applicability of the OT to establish sharp Sobolev inequalities on non-euclidean spaces, they wrote 
 (see \cite[p.\ 309]{CE-N-Villani}): 
 %\enquote*{A}
 \enquote{(...) \textit{We do not know whether our methods would still be as efficient in a Riemannian setting}.}\
  Closely related to this question, Villani \cite[Chapter 21]{Villani} has formulated several open problems concerning the possible proof of sharp geometric/functional inequalities on  non-euclidean structures using only  OT theory.

Motivated by Villani's questions, 
  Cavalletti and Mondino \cite{C-M-inv, C-M-2} proved the sharp L\'evy-Gromov isoperimetric inequality as well as sharp Sobolev and logarithmic-Sobolev inequalities on compact ${\sf CD}(K,N)$ metric measure spaces with $K>0$   via OT-arguments; see Lott and Villani \cite{LV} and Sturm \cite{Sturm} for the theory of these metric measure spaces.
%  By  applying $1$-dimensional localization, Cavalletti and Mondino \cite{C-M-2} established sharp Sobolev and logarithmic-Sobolev inequalities for ${\sf CD}(K,N)$ spaces with finite diameter. 
%  
%  Cavalletti and Mondino \cite{C-M-inv} provided the first OT-based proof of the sharp L\'evy-Gromov isoperimetric inequality on ${\sf CD}(k,N)$ metric measure spaces with $k>0$; for basic properties of these metric measure spaces, see Lott and Villani \cite{LV} and Sturm \cite{Sturm}.
   The  approach in \cite{C-M-inv, C-M-2} is based on Klartag's disintegration technique, see \cite{Klartag}, which is an $L^1$-OT-argument, whose origin  goes back to the localization paradigm from convex geometry developed by  Gromov and Milman  \cite{G-M}, Kannan, Lov\'asz and Simonovits \cite{KLS} and Lov\'asz and  Simonovits \cite{LS}. The latter technique -- similarly to symmetrization -- reduces multidimensional problems to  $1$-dimensional model objects.

   Very recently, by using the $L^2$-OT-based distorted Brunn-Minkowski inequality, Balogh and Krist\'aly \cite{BK} stated the sharp isoperimetric inequality on  non-compact  ${\sf CD}(0,N)$ spaces, extending previous contributions by Agostiniani,  Fogagnolo and  Mazzieri \cite{AFM}, Brendle \cite{Brendle} and Fogagnolo and  Mazzieri \cite{FM}. As a consequence, following the technique of Aubin \cite{Aubin} and Talenti \cite{Talenti},  the latter sharp  isoperimetric inequality combined with suitable  symmetrization  provided the sharp $L^p$-Sobolev inequality $(p>1)$ on Riemannian manifolds with ${\sf Ric}\geq 0$, see \cite[Theorem 1.2]{BK}. 
 
 Our first goal is to provide a \textit{symmetrization-free, fully $L^2$-OT-based} proof of the  sharp $L^p$-Sobolev inequality, $p\geq 1,$ answering in this way the question of Cordero-Erausquin,  Nazaret and Villani \cite{CE-N-Villani} on noncompact, complete Riemannian manifolds with ${\sf Ric}\geq 0$. 
 
 To be more precise, let $(M,g)$ be a noncompact, complete  $n(\geq 2)$-dimensional Riemannian manifold with ${\sf Ric}\geq 0$, endowed with its volume form ${\rm d}v_g$ and distance function $d_g$.  Let
 $${\sf AVR}_g=\lim_{r\to \infty}\frac{{\rm Vol}_g(B_x(r))}{\omega_nr^n}
 $$
 be the 
 {\it asymptotic volume ratio} on $(M,g)$;
 due to Bishop-Gromov  comparison, ${\sf AVR}_g$
 is well-defined, independent of the choice of $x\in M$, and ${\sf AVR}_g\in [0,1]$. Here,  $B_x(r)$ is the metric ball with center at $x\in M$ and radius $r>0$, ${\rm Vol}_g(S)$ being the volume of the set $S\subset M$. For further use, $\nabla_g$ stands for the Riemannian gradient on $(M,g)$, $C_0^\infty(M)$ is the space of  compactly supported smooth functions on $M$, while for a function $f:M\to \mathbb R$, the number $|\nabla_g f(x)|$ is  the norm of $\nabla_g f(x)\in T_xM$ in the Riemannian metric over the tangent space $T_xM$, $x\in M.$ 
 
% constant ${\sf AVR}_g\in (0,1]$,

% Moreover, they raised the question on the possible extension of the OT-approach to prove sharp Sobolev inequalities on Riemannian manifolds. More precisely, they wrot
 
%\begin{quotation}
%	"\textit{Considerable effort has been spent recently on the problem of optimal Sobolev
%		inequalities on Riemannian manifolds, see the survey \cite{Druet-Hebey} and references therein. In
%		the present work however, we shall concentrate on the situation where the problem is
%		set on $\mathbb R^n.$ We do not know whether our methods would still be as efficient in a
%		Riemannian setting. Note however that nonsharp Sobolev Riemannian inequalities
%		can easily be derived by mass transportation techniques, as shown in \cite{CE-thesis}.}"
%\end{quotation}
%This question is closely related to the AB-program initiated by Aubin \cite{Aubin} to determine the best constants in various Sobolev inequalities on Riemannian manifolds.  	

%As we already pointed out, our purpose is to  answer the aforementioned question for the class of Riemannian manifolds with ${\sf Ric}\geq 0$. 

% namely, \textit{by using exclusively OT-theory} (and some elementary inequalities, as H\"older, AM-GM) we are going to prove  \textit{sharp $L^p$-Sobolev inequalities on Riemannian manifolds with ${\sf Ric}\geq 0$}. 

	The \textit{sharp $L^p$-Sobolev inequality} $(p\geq 1)$, which we prove using  OT-arguments, reads as follows:

\begin{theorem}\label{main-thm}  Let $(M,g)$ be a noncompact, complete  $n$-dimensional Riemannian manifold $(n\geq 2)$ with ${\sf Ric}\geq 0$ having Euclidean volume growth, i.e., $0< {\sf AVR}_g \leq 1$.  If $p\in [1,n)$, then  
	\begin{equation}\label{egyenlet-0}
		\displaystyle  { \left(\int_{M} |f|^{p^\star}{\rm d}v_g\right)^{1/p^\star}\leq
			{\sf AT}(n,p)\, {\sf AVR}_g^{-\frac{1}{n}}
			\left(\int_{M} |\nabla_g f|^p {\rm d}v_g\right)^{1/p}},\ \ \forall  f\in C_0^\infty(M).
	\end{equation}
	Moreover, the constant $ %{\sf S}_g=
	{\sf AT}(n,p)\, {\sf AVR}_g^{-\frac{1}{n}}$ in \eqref{egyenlet-0} is   sharp. 
\end{theorem}

We notice that in the setting of Theorem \ref{main-thm}, the Euclidean volume growth requirement, i.e., ${\sf AVR}_g>0$, is nothing but a non-collapsing assumption that guarantees the validity of the underlying Sobolev embedding, see e.g.\ Coulhon and Saloff-Coste \cite{C-SC} and Hebey \cite[Theorem 3.2]{Hebey}, and is also necessarily  for a Sobolev inequality of the type \eqref{egyenlet-0} to hold, see e.g.\ Ledoux \cite{Ledoux} and Xia \cite{Xia}.  

In our approach
we use deep $L^2$-OT-arguments developed by McCann \cite{McCann} and Cordero-Erausquin,  McCann and Schmuckenschl\"{a}ger \cite{CEMS}; the key tools are as follows: 
\begin{itemize}
	\item \textit{Determinant-trace inequality} (see Proposition \ref{Wang-Zhang-proposition}). This inequality is  derived from the representation  of the Jacobian of the optimal transport map between the  marginal measures, 
	together with certain comparison arguments coming from the fact that ${\sf Ric}\geq 0$.
	\item \textit{Monge-Amp\`ere equation}. Properties of the $c=d_g^2/2$-concave function arising from the  optimal transport map  are  combined with the Monge-Amp\`ere equation and suitable change of variables.  
	\item \textit{Asymptotic properties of Talentian bubbles}  (see Proposition \ref{AT-estimates-proposition}).  This argument produces the appearance of ${\sf AVR}_g$ jointly with the Aubin--Talenti-constant  ${\sf AT}(n,p)$ in  \eqref{egyenlet-0}, where the property  ${\sf Ric}\geq 0$  plays again an indispensable role. 
\end{itemize}
   In addition, the asymptotic property of Talentian bubbles -- which incidentally provides the sharp Sobolev constant in \eqref{egyenlet-0} -- 
 yields \textit{per se} an elegant, alternative proof  of the non-collapsing result of do Carmo and Xia \cite{doCarmo-Xia} and other related results, see \S \ref{section-5}.

Our paper focuses on the sharpness, the equality case in \eqref{egyenlet-0} not being discussed. The reason is that our approach contains limiting arguments that require to work  with functions belonging  to $ C_0^\infty(M)$, even though  the expected 
\enquote*{candidates}
 for the extremals in \eqref{egyenlet-0} are contained in a larger Sobolev space (as in the Euclidean case). 
% -- although both \eqref{egyenlet-0} and \eqref{LSI} can be extended in the usual way from $ C_0^\infty(M)$  to suitable Sobolev spaces --  and the latter space does not contain the expected 'candidates' for the extremals in \eqref{egyenlet-0} and \eqref{LSI}.
However,  if equality is achieved by a (positive) extremal in \eqref{egyenlet-0},  this function solves an elliptic PDE involving the critical exponent $p^\star$. The latter class of problems has been widely studied  recently by Fogagnolo,  Malchiodi and Mazzieri \cite{olaszok-1} as well as by Catino, Monticelli and Roncoroni \cite{CM, CMR} (under some mild technical assumptions/restrictions), stating that a Riemannian manifold $(M,g)$ with  ${\sf Ric}\geq 0$ supports the solvability of such elliptic PDEs if and only if %the Riemannian manifold 
$(M,g)$ is isometric to the Euclidean space; in particular, the extremal functions are precisely the Talentian bubbles $(f_\lambda)_{\lambda>0}$ from \eqref{talentian-bubble-eukl}. Moreover, in a very recent paper, Nobili and Violo \cite[Theorem 5.3]{NV-adv} characterized the equality case in  \eqref{egyenlet-0} for  $p=2$ even on non-smooth ${\sf RCD}(0,N)$  spaces.

To further emphasize the efficiency of our OT-arguments, we prove the \textit{sharp $L^p$-logarithmic Sobolev inequality} $(p\geq 1)$ on Riemannian manifolds with  ${\sf Ric}\geq 0$. To state this result, let $p>1$, $n\geq 2$ and $$	
{\sf L}({n,p})={\small
	\frac{p}{n}\left(\frac{p-1}{e}\right)^{p-1}\left(\omega_n{\Gamma\left(\frac{n}{p'}+1\right)}\right)^{-\frac{p}{n}}}. 
$$
For the limit case $p=1$, let ${\sf L}(n,1)=\lim_{p\searrow 1}{\sf L}(n,p)=n^{-1}\omega_n^{-\frac{1}{n}}={\sf AT}(n,1).$ 

Our second main result can be stated as follows: 

%approach turns out to be efficient also to provide a direct, OT-based proof for the sharp log-Sobolev inequality; namely, we have:

	\begin{theorem}\label{log-Sobolev-main} Let $(M,g)$ be an $n$-dimensional  Riemannian manifold as in Theorem \ref{main-thm} and $p\geq 1$. Then  for every $ f\in C_0^\infty(M)$ 
		 with $\displaystyle\int_M |f|^p {\rm d}v_g=1$, one has that 
	\begin{equation}\label{LSI}
		\int_{M}|f|^p\log |f|^p{\rm d}v_g\leq \frac{n}{p}\log\left({\sf L}({n,p}){\sf AVR}_ g^{-\frac{p}{n}}\int_M |\nabla_g f|^p{\rm d}v_g \right). 
	\end{equation}
	Moreover,  the constant  ${\sf L}({n,p}){\sf AVR}_g^{-\frac{p}{n}}$ in \eqref{LSI} is sharp.  
\end{theorem}

Similarly as in Theorem \ref{main-thm}, the proof of Theorem \ref{log-Sobolev-main} is based on the \textit{determinant-trace inequality} and on fine asymptotic properties of \textit{Gaussian bubbles} on Riemannian manifolds with  ${\sf Ric}\geq 0$ (see Proposition \ref{L-estimates-proposition}), the latter argument providing again the sharp constant ${\sf L}({n,p}){\sf AVR}_g^{-\frac{p}{n}}$ in \eqref{LSI}.

We notice that the sharp $L^p$-logarithmic Sobolev inequality in the Euclidean space $\mathbb R^n$ has been established first by Weissler \cite{Weissler} for $p = 2,$
Del Pino and Dolbeault \cite{delPino} for $p\in (1,n)$, and Gentil \cite{Gentil} for general $p>1$. Furthermore, by using the symmetrization/rearrangement-type argument from Nobili and Violo \cite{NV}, the sharp $L^p$-logarithmic Sobolev inequality with $p>1$ has been extended very recently by Balogh, Krist\'aly and Tripaldi \cite{BKT} to  ${\sf CD}(0,N)$ spaces.

%  very recent studies, based on elliptic PDE theory containing the critical exponent,    proved that under certain technical restrictions, 

We conclude the presentation of our results by a \textit{Gaussian logarithmic-Sobolev inequality}, which turns out to be a direct consequence of Theorem \ref{log-Sobolev-main} for $p=2$. To state it, let
 $(M,g)$ be an $n$-dimensional Riemannian manifold and  
 $V\in C^2(M)$ verifying 
 \begin{equation}\label{assumption-Cv}
 	V-\frac{|\nabla_g V|^2}{2}+\Delta_g V-n\leq C_V\ \ {\rm a.e.\ on} \ M,
 \end{equation}
 for some constant $C_V\in \mathbb R$, where  $\Delta_g$ stands for the Laplace-Beltrami operator on $(M,g)$.

%  Let ${\rm d}\gamma_V=G_V^{-1}e^{-V}{\rm d}v_g$ be a probability measure on the Riemannian manifold  with  and  satisfying 

\begin{theorem}\label{Gaussian-log-Sobolev-main} Let $(M,g)$ be an $n$-dimensional  Riemannian manifold as in Theorem \ref{main-thm}. Assume that  $V\in C^2(M)$ satisfies  \eqref{assumption-Cv} and $G_V=\int_M e^{-V}{\rm d}v_g<\infty$,  ${\rm d}\gamma_V=G_V^{-1}e^{-V}{\rm d}v_g$ being a probability measure on $(M,g)$. 
	Then  for every $ h\in C_0^\infty(M)$ 
	with $\displaystyle\int_M h^2 {\rm d}\gamma_V=1$, one has that 
	\begin{equation}\label{Gaussian-LSI}
		\int_{M}h^2\log h^2{\rm d}\gamma_V\leq 2\int_M |\nabla_g h|^2{\rm d}\gamma_V+\log\left(\frac{G_Ve^{C_V}}{(2\pi)^\frac{n}{2} {\sf AVR}_g}\right).
	\end{equation} 
\end{theorem}
A typical example verifying the assumptions is  $V(x)=\frac{1}{2}d_g^2(x_0,x)$ for some  $x_0\in M$ (see Remark \ref{remark-particular} for details), whenever \eqref{Gaussian-LSI} reduces to the more familiar, \enquote*{dimension-free} Gaussian logarithmic-Sobolev inequality
	\begin{equation}\label{Gaussian-LSI-particular}
	\int_{M}h^2\log h^2{\rm d}\gamma_V\leq 2\int_M |\nabla_g h|^2{\rm d}\gamma_V-\log{\sf AVR}_g.
\end{equation} 
%\textbf{Semilinear elliptic equations on manifolds with nonnegative Ricci curvature}

The paper is organized as follows. In \S \ref{section-2} we prepare the proofs of our main results. First, in  \S \ref{det-ineq-section} we prove the key determinant-trace inequality (see Proposition \ref{Wang-Zhang-proposition}) where we also recall some basic notions and results from the theory of OT. Then, in \S\ref{subsection-bubbles} we establish the fine asymptotic behaviors of the Talentian and Gaussian bubbles on Riemannian manifolds with ${\sf Ric}\geq 0$ (see Propositions \ref{AT-estimates-proposition} \& \ref{L-estimates-proposition}).  In \S\ref{section-3}	we prove the sharp $L^p$-Sobolev inequality (see Theorem \ref{main-thm}) for both $p>1$  and $p=1$; although the  case $p=1$ follows by limiting from  $p>1$, we provide an independent proof. In \S\ref{section-4}  the sharp $L^p$-logarithmic Sobolev inequality (see Theorem \ref{log-Sobolev-main}) is proven for both $p>1$  and $p=1$, followed by the Gaussian logarithmic-Sobolev inequality (Theorem \ref{Gaussian-log-Sobolev-main}). Finally, \S \ref{section-5} 
is devoted to an elementary proof of the main result of do Carmo and Xia \cite{doCarmo-Xia}.

\section{Preliminaries}\label{section-2}

\subsection{Determinant-trace inequality}\label{det-ineq-section} The purpose of this subsection is to prove a  determinant-trace inequality based on OT that will be of crucial importance in the proof of our main results. 

Let  $(M,g)$ be a  complete $n$-dimensional Riemannian manifold $(n\geq 2)$ with the induced distance function $d_g:M\times M\to [0,\infty)$ and canonical volume form ${\rm d}v_g$.
%We first recall some notions from the OT-theory  needed to carry out our proofs; for details, see Cordero-Erausquin, McCann and  Schmuckenschl\"{a}ger \cite{CEMS} and McCann \cite{McCann}. 
%and
For any $t\in [0,1]$ and $x,y\in M$, let 
$$Z_t(x,y)=\{z\in M: d_g(x,z)=td_g(x,y)\ \ {\rm and}\ \ d_g(z,y)=(1-t)d_g(x,y)\}$$
be the set of  \textit{$t$-intermediate points} between $x$ and $y.$
Since
$(M,d_g)$ is complete, we clearly have that $Z_t(x,y)\neq \emptyset$.
Let $B_x(r)=\{y\in
M:d_g(x,y)<r\}$ be the geodesic ball with center $x\in M$ and radius  $r>0$. Fix $t\in (0,1).$ 
According to Cordero-Erausquin, McCann and Schmuckenschl\"ager
\cite{CEMS},  the \textit{volume distortion} in
$(M, g)$ is defined by
$$
	v_t(x,y) = \lim\limits_{r \to 0}\frac{{\rm Vol}_g\left(Z_t(x, B_y( r))\right)}{{\rm Vol}_g\left(B_y(tr)\right)},
$$
where ${\rm Vol}_g(S)$ denotes the volume of the set $S\subset M.$
In particular, when ${\sf Ric}\geq 0$, it turns out that for every $t\in (0,1)$ and  $x,y\in M$ with $y\notin {\sf cut}(x)$,  
\begin{equation}\label{v-s-nagyobb-mint-1}
	v_t(x,y)\geq 1,
\end{equation}
see \cite[Corollary 2.2]{CEMS}, where ${\sf cut}(x)\subset M$ denotes the cut-locus of $x.$

Let $\Omega\subset M$ be an open bounded set, and let  ${\rm d}\mu(x)= F(x) {\rm d}v_g(x)$ and ${\rm d}\nu(y)= G(y){\rm d}
v_g(y)$ be two absolutely continuous probability measures on $(M,g)$ with ${\rm supp}\mu=\overline \Omega$. According to McCann \cite{McCann}, there exists a unique  optimal transport
map $T:\overline \Omega\to M$  with $\nu=T_\#\mu$ (i.e., $\nu$ is the push-forward of $\mu$ through $T$), having  the form $$T(x)=\exp_x(-\nabla_g u(x))\ \ {\rm for\ a.e.}\ \  x\in \operatorname{supp} \mu,$$ where $u:\overline \Omega\to \mathbb R$  is a $c=d_g^2/2$-concave function and $\nabla_g$ denotes the Riemannian gradient on $(M,g)$.  We recall that $u:\overline \Omega\to \mathbb R$  is a $c=d_g^2/2$-concave function if there exists a function $\eta:Y\to \mathbb R\cup \{-\infty\}$ with $\emptyset \neq \operatorname{supp} \nu\subseteq Y\subset M$ such that
$$u(x)=\inf_{y\in Y}\left(\frac{d_g^2(x,y)}{2}-\eta(y)\right),\ x\in M.$$
If $u:\overline \Omega\to \mathbb R$ is $c=d_g^2/2$-concave, then it is semi-concave on $\Omega$, and hence admits a Hessian a.e.\ in $\Omega$, see \cite[Proposition 3.14]{CEMS}. Furthermore, the push-forward relation $T{_\sharp} \mu=\nu$ implies the \textit{Monge-Amp\`ere equation}, i.e.,  
$$
	F(x)= G(T(x)){\rm det}DT(x)\ {\rm for \ \mu}-{\rm a.e.}\ x\in \Omega.
$$
%Let $s\in (0,1)$ and $T_s:M\to M$ be the $s$-interpolant optimal transport map defined by $$T_s(x)=\exp_x(-s\nabla_g u(x))\ {\rm for \ a.e.}\ x\in \overline\Omega.$$ One of the most important relations that we shall use throughout this section is the following
%Jacobian determinant inequality  \cite[Lemma 6.1]{CEMS} that holds for a.e. $x\in \Omega$ and reads as follows:
%\begin{equation}\label{Jacobian-inequality}
%	\left( \det DT_{s} (x) \right)^{\frac{1}{n}} \geq (1-s) v_{1-s} (T(x), x))^{\frac{1}{n}}  + s v_s(x, T(x)) ^{\frac{1}{n}}  \left( \det DT(x)\right)^{\frac{1}{n}}.
%\end{equation}

%In order to formulate our first statement, we introduce the functions $\mathcal H(s)=s\coth(s)$ and $\mathcal S(s)=\sinh(s)/s$, $s>0$. 
A key result in our arguments is the following determinant-trace inequality: 
\begin{proposition}\label{Wang-Zhang-proposition}
	Let $(M,g)$ be a noncompact, complete  $n$-dimensional Riemannian mani\-fold $(n\geq 2)$ with ${\sf Ric}\geq 0 $, $\Omega\subset M$ be an open bounded set and $T:\overline \Omega\to M$ be the unique optimal transport map between two probability measures $\mu$ and $\nu$ on $M$ with ${\rm supp}\mu=\overline\Omega$, given by  $T(x)=\exp_x(-\nabla_g u(x))$ for a.e. $x\in \Omega$ for some $c=d_g^2/2$-concave function $u:\overline\Omega\to \mathbb R$. Then for $\mu$-a.e.\ $x\in \Omega$ one has 
		\begin{equation} \label{eq-WZ-1} 
			({\rm det}DT(x))^\frac{1}{n} \leq 1-\frac{\Delta_g u(x)}{n}.
		\end{equation} 
%		\item[(ii)]
%		If $\kappa >0,$ then for a.e. $x\in \Omega$ it holds 
%		\begin{equation} \label{eq-WZ-2} 
%			{\rm det}DT(x) \leq \mathcal S^{n}\left(\sqrt{\frac{\kappa}{n}}|\nabla_gu(x)|\right) \left(\mathcal H\left(\sqrt{\frac{\kappa}{n}}|\nabla_gu(x)|\right)-\frac{\Delta_g u(x)}{n}\right)^n. 
%		\end{equation} 
%		
%	\end{enumerate} 
\end{proposition}

{\it Proof.} 
Due to Cordero-Erausquin, McCann and Schmuckenschl\"ager
\cite[Theorem 4.2]{CEMS} we have the representation of  Jacobian of the optimal transport map $T$, given by
\begin{equation}\label{Jacobian-repr}
	DT(x)=Y_x(1)H(x)\ \ {\rm for}\ \ \mu-{\rm  a.e.}\   x\in \Omega,
\end{equation}
where $$Y_x(1)=d(\exp_x)_{(-\nabla_g u(x))}\ \ {\rm and}\ \ H(x)= {\rm Hess}_g \left[\frac{1}{2}d^2_y(x) - u(x)\right]$$ with $d_y(x)=d_g(y,x)$ and $y=T(x)$. Let us recall that in relation \eqref{Jacobian-repr} it is implicitly contained the fact that $T(x) = y\notin {\sf cut}\{x\}$
for $\mu$-a.e.\ $x\in \Omega$, which in particular implies the fact that $x \mapsto d_g^2(x, y)$ is smooth. In the sequel we shall always work with such points $x$ that form a full measure set. 

 Since ${\sf Ric}\geq 0$, the Laplace comparison principle implies that  $\Delta_g d^2_y(x) \leq 2n $ for $y=T(x)$, see e.g.\ Petersen \cite[Lemma 7.1.2]{Petersen}.   Moreover,  by 
\cite[Proposition 4.1]{CEMS}, it turns out that for $\mu$-a.e.\ $x\in \Omega$ the matrix $H(x)$ is nonnegative definite and symmetric. Thus, by the AM-GM inequality  we obtain that
\begin{eqnarray*}
	 ({\rm det}H(x))^\frac{1}{n}   &\leq&  \frac{1}{n} {\rm trace} H(x) =   \frac{1}{n}\left\{ \frac{1}{2}\Delta_g d^2_y(x) - \Delta_ g u(x) \right\}\\& \leq& 1-\frac{\Delta_g u(x)}{n} .
\end{eqnarray*}
	On the other hand,  if we consider $Y_x(t)=d(\exp_x)_{(-t\nabla_g u(x))}$ for $t\in (0,1),$  by  \cite[Lemma 2.1]{CEMS} and  \eqref{v-s-nagyobb-mint-1}, we have that
$$\frac{{\rm det}Y_x(t)}{{\rm det}Y_x(1)}=v_t(x,T(x))\geq 1\ \ {\rm for}\ \ \mu-{\rm  a.e.}\   x\in \Omega.$$
In particular, it follows that  for $\mu$-a.e.\ $x\in \Omega  $ and every $t\in (0,1)$ one has ${\rm det}Y_x(1)\leq {\rm det}Y_x(t).$
Since $d(\exp_x)_0={\rm Id}$, see e.g.\ do Carmo \cite[Proposition 2.9]{doCarmo}, once we take the limit $t\to 0^+$ in the latter relation, we obtain  that ${\rm det}Y_x(1)\leq 1,$
which concludes the proof of \eqref{eq-WZ-1}. 
\hfill $\square$

\begin{remark}\rm 
	Based on the Alexandroff--Bakelman--Pucci-theory, a similar inequality to \eqref{eq-WZ-1} has been proved by Wang and Zhang \cite[Theorem 1.2]{WZ} for functions $u\in  C^2(\Omega)$ verifying a Savin-type \textit{contact-property}. In fact, a closer inspection of the definitions shows that the contact-property, up to a sign-convention, is equivalent to the $c$-concavity. 
\end{remark}

The determinant-trace inequality \eqref{eq-WZ-1} contains the Laplacian $\Delta_g$,  
%of the $c$-concave function $u$ in the sense of Alexandrov, 
which is the absolutely continuous part of the distributional Laplacian $\Delta_{g,\mathcal D'}$.  We now provide  a technical result connecting these operators, which will be important in the application of the divergence theorem. 

\begin{proposition}\label{Laplace-singular}
	Let $(M,g)$ be a complete Riemannian manifold,  $\Omega\subset M$ be an open set and $u:\Omega\to \mathbb R$ be a $d_g^2/2$-concave function on $\Omega$. Then 
	the singular part  $\Delta_g^{\rm s} (-u)$ of the distributional Laplacian $\Delta_{g,\mathcal D'}(-u)$ is a positive Radon measure on $\Omega$. In particular, $-\Delta_g u\leq -\Delta_{g,\mathcal D'} u$. 
\end{proposition}

{\it Proof.} Since $u:\Omega\to \mathbb R$ is a $d_g^2/2$-concave function on $\Omega$, it is also semi-concave, see Cordero-Erausquin, McCann and Schmuckenschl\"ager
\cite[Proposition 3.14]{CEMS}. In particular,  to every point $x \in \Omega$ there exist $C>0$ and enough small $r_x>0$  
such that $C\frac{d_g^2(x,\cdot)}{2}-u$ is geodesically convex on $B_x(r_x)$ (i.e., it is convex in the usual sense along any geodesic segment belonging to $B_x(r_x)$). By adding to this function a multiple of 
$d_g^2(x,\cdot)$ we can assume without loss of generality that the previous function is even geodesically strictly convex. This property  implies that %the singular part of the distributional derivative 
${\rm div}^{\rm s} \left(\nabla_g\left(C\frac{d_g^2(x,\cdot)}{2}-u\right)\right)$ is a positive Radon measure on $B_x(r_x)$. Since $d_g^2(x,\cdot)$ is smooth on $B_x(r_x)$, it follows that 
$${\rm div}^{\rm s} \left(\nabla_g\left(C\frac{d_g^2(x,\cdot)}{2}-u\right)\right) = {\rm div}^{\rm s} \left(\nabla_g\left(-u\right)\right)=\Delta_g^{\rm s} (-u)$$ is a positive Radon measure on $B_x(r_x)$. Using the local property of distributional derivatives (applying an appropriate partition of unity on the set $\Omega$), we obtain that the singular part of the distributional Laplacian is a positive Radon measure on the whole set $\Omega$. In particular, we have that
$$-\Delta_{g,\mathcal D'} u=\Delta_{g,\mathcal D'}(-u)=-\Delta_{g} u+\Delta_g^{\rm s} (-u)\geq -\Delta_{g} u,$$
which concludes the proof. 
\hfill $\square$

%\textcolor{red}{Would it be good to give a more detailed proof of the second part and state it in its sharper form?}

	\subsection{Bubbles on Riemannian manifolds with ${\sf Ric}\geq 0$} \label{subsection-bubbles}
	Let $(M,g)$ be a complete  $n$-dimensional Riemannian manifold $(n\geq 2)$ with ${\sf Ric}\geq 0 $, and $x_0\in M$ be fixed.  
	Given $\lambda,s>0$ and $p\in (1,n),$ we consider the \textit{Talenti bubble} $x\mapsto (\lambda+d_g^{p'}(x_0,x))^{-s}$ and its associated integral function 
	\begin{equation}\label{H-function-definition}
			\mathcal H(\lambda,s)=\int_M(\lambda+d_g^{p'}(x_0,x))^{-s}{\rm d}v_g(x).
	\end{equation}

	In the sequel, we collect the basic properties of $\mathcal H$:

	\begin{proposition}\label{AT-estimates-proposition}
		Let $(M,g)$ be a noncompact, complete  $n$-dimensional Riemannian mani\-fold $(n\geq 2)$ with ${\sf Ric}\geq 0 $. If $p\in (1,n)$ and $s>\frac{n}{p'},$ then 
		\begin{itemize}
			\item[(i)] $0<\mathcal H(\lambda,s)<\infty$ for every $\lambda>0;$
			\item[(ii)]  the following  asymptotic property holds$:$ 
				\begin{equation} \label{eq-AT-1} 
				\lim_{\lambda\to \infty} \lambda^{s-\frac{n}{p'}} \mathcal H(\lambda,s)= \omega_n{\sf AVR}_g\frac{\Gamma\left(\frac{n}{p'}+1\right)\Gamma\left(s-\frac{n}{p'}\right)}{\Gamma(s)}.
			\end{equation} 
		\end{itemize}
\end{proposition}		
	{\it Proof.} 
	(i) By the layer cake representation and a change of variables one easily obtains that 
\begin{equation}\label{H-estimate}
	\mathcal H(\lambda,s)=\int_M(\lambda+d_g^{p'}(x_0,x))^{-s}{\rm d}v_g(x)=sp'\int_0^\infty{\rm Vol}_g(B_{x_0}(\rho))(\lambda+\rho^{p'})^{-s-1}\rho^{p'-1}{\rm d}\rho.
\end{equation}
	Since ${\sf Ric}\geq 0 $, the Bishop-Gromov volume comparison principle implies that ${\rm Vol}_g(B_{x_0}(\rho))\leq \omega_n\rho^n$ for every $\rho>0$; thus,  by our assumption $s>\frac{n}{p'}$ it follows that
		$$\mathcal H(\lambda,s)\leq sp'\omega_n\int_0^\infty(\lambda+\rho^{p'})^{-s-1}\rho^{n+p'-1}{\rm d}\rho<\infty.$$
		
		(ii) A change of variables of the form $\rho=\lambda^{\frac{1}{p'}}r$ in \eqref{H-estimate} gives that 
		\begin{eqnarray*}
				\mathcal H(\lambda,s)&=&\lambda^{-s}sp'\int_0^\infty{\rm Vol}_g\left(B_{x_0}(\lambda^{\frac{1}{p'}}r)\right)(1+r^{p'})^{-s-1}r^{p'-1}{\rm d}r\\&=&\omega_n\lambda^{-s+\frac{n}{p'}}sp'\int_0^\infty\frac{{\rm Vol}_g\left(B_{x_0}(\lambda^{\frac{1}{p'}}r)\right)}{\omega_n\left(\lambda^{\frac{1}{p'}}r\right)^n} (1+r^{p'})^{-s-1}r^{n+p'-1}{\rm d}r. 
		\end{eqnarray*}
		Taking the limit in  $\lambda^{s-\frac{n}{p'}} \mathcal H(\lambda,s)$ as $\lambda\to \infty$, the monotone convergence theorem (as $\rho\mapsto \frac{{\rm Vol}_g(B_{x_0}(\rho))}{\rho^n}$ is non-increasing by the Bishop-Gromov comparison principle) and the definition of ${\sf AVR}_g$ imply
		$$\lim_{\lambda\to \infty} \lambda^{s-\frac{n}{p'}} \mathcal H(\lambda,s)=\omega_n{\sf AVR}_gsp' \int_0^\infty (1+r^{p'})^{-s-1}r^{n+p'-1}{\rm d}r.$$ It remains to use basic properties of the Gamma function $\Gamma$ to conclude the proof of \eqref{eq-AT-1}. 
	\hfill $\square$\\
		
	In a similar way as before, we consider the \textit{Gaussian bubble} $x\mapsto e^{-\lambda d_g^{p'}(x_0,x)}$ and the functions 
	\begin{equation}\label{LL-function-definition}
		\mathcal L_1(\lambda)=\int_Me^{-\lambda d_g^{p'}(x_0,x)}{\rm d}v_g(x)\ \ {\rm and}\ \ \mathcal L_2(\lambda)=\int_Me^{-\lambda d_g^{p'}(x_0,x)}d_g^{p'}(x_0,x){\rm d}v_g(x).
	\end{equation}
	
		\begin{proposition}\label{L-estimates-proposition}
		Let $(M,g)$ be a noncompact, complete  $n$-dimensional Riemannian mani\-fold $(n\geq 2)$ with ${\sf Ric}\geq 0 $. If $p>1$  then 
		\begin{itemize}
			\item[(i)] $0<\mathcal L_i(\lambda)<\infty$ for every $\lambda>0$ and $i\in \{1,2\},$ and $\mathcal L_2 =-\mathcal L_1';$
			\item[(ii)]  the following  asymptotic properties hold$:$ 
			$$
				\lim_{\lambda\to 0} \lambda^{\frac{n}{p'}} \mathcal L_1(\lambda)= \omega_n{\sf AVR}_g{\Gamma\left(\frac{n}{p'}+1\right)}\ \ {and}\ \ \lim_{\lambda\to 0} \lambda^{\frac{n}{p'}+1} \mathcal L_2(\lambda)= \omega_n{\sf AVR}_g\frac{n}{p'}{\Gamma\left(\frac{n}{p'}+1\right)}.
			$$
		\end{itemize}
	\end{proposition}		
	{\it Proof.} We only deal with $\mathcal L_1$, the proof being similar also for $\mathcal L_2$. 
	
	(i) By using the layer cake representation, it follows that $$\mathcal L_1(\lambda)=\int_Me^{-\lambda d_g^{p'}(x_0,x)}{\rm d}v_g(x)=\lambda p'\int_0^\infty{\rm Vol}_g(B_{x_0}(\rho))e^{-\lambda \rho^{p'}}\rho^{p'-1}{\rm d}\rho.$$
	In particular, due to the Bishop-Gromov volume comparison principle, $\mathcal L_1$ is well defined. 
	
(ii) By the change of variables $\rho=\lambda^{-\frac{1}{p'}}r^{\frac{1}{p'}}$, it follows that $$\mathcal L_1(\lambda)=\int_0^\infty{\rm Vol}_g\left(B_{x_0}\left(\lambda^{-\frac{1}{p'}}r^{\frac{1}{p'}}\right)\right)e^{-r}{\rm d}r.$$ Therefore, by  the monotone convergence theorem we obtain 
$$	\lim_{\lambda\to 0} \lambda^{\frac{n}{p'}} \mathcal L_1(\lambda)= \omega_n{\sf AVR}_g\int_0^\infty e^{-r}r^\frac{n}{p'} {\rm d}r= \omega_n{\sf AVR}_g{\Gamma\left(\frac{n}{p'}+1\right)},$$
which is the desired relation. 
		\hfill $\square$
	
	\section{Proof of the sharp $L^p$-Sobolev inequality (Theorem \ref{main-thm})}\label{section-3}
\subsection{The case $p>1$} Let us fix arbitrarily 	$ f\in C_0^\infty(M)$; in order to prove \eqref{egyenlet-0}, we may assume without loss of generality that $f$ is nonnegative and for simplicity, 
\begin{equation}\label{normalized-1}
	\int_{M} f^{p^\star}{\rm d}v_g=1.
\end{equation}
Let $\Omega=\{x\in M:f(x)>0\}$; since $ f\in C_0^\infty(M)$, then $\overline \Omega$ is compact. 
 
 For further use, for every $k\in \mathbb N$ we consider the truncation function $P_k:[0,\infty)\to \mathbb R$ defined by
 \begin{equation}\label{truncation-k}
 	P_k(s)=\max(0,\min(0,k-s)+1).
 \end{equation}
It is clear that $P_k\leq P_{k+1}$ and the support of $P_k$ is $[0,k+1]$. 

Let $x_0\in \Omega.$ Fix $\lambda>0$, $k\in \mathbb N$ arbitrarily and consider the truncated Talentian bubble   $G_{\lambda,k}:M\to \mathbb R$ given by
\begin{equation}\label{G-extrem-approx}
	G_{\lambda,k}(x)=P_k(d_g(x_0,x))\left(\lambda+d_g^{p'}(x_0,x)\right)^{-n}.
\end{equation}
By construction, the support of $G_{\lambda,k}$ is  $\overline {B_{x_0}(k+1)}$. Let $$\mathcal I_{\lambda,k}=\int_{M}G_{\lambda,k}(y){\rm d}v_g(y).$$
	It is clear that $0<\mathcal I_{\lambda,k}<\infty$ for every $\lambda>0$ and $k\in \mathbb N$.   
	
Let us consider the compactly supported probability measures ${\rm d}\mu(x)= f(x)^{p^\star} {\rm d}v_g(x)$ and ${\rm d}\nu(y)= \frac{G_{\lambda,k}(y)}{\mathcal I_{\lambda,k}}{\rm d}
	v_g(y)$  on $(M,g)$.   Due to McCann \cite{McCann}, there exists a unique  optimal transport
	map $T:\overline \Omega\to \overline {B_{x_0}(k+1)}\subset M$  pushing $\mu$ forward to $\nu$ and  
\begin{equation}\label{Optim-map}
	T(x)=\exp_x(-\nabla_g u(x))\ \ {\rm for\ a.e.}\ \  x\in \overline \Omega,
\end{equation}	
	 where $u:\overline \Omega\to \mathbb R$  is a $c=d_g^2/2$-concave function. Clearly, both $T$ and $u$ depend on $\lambda>0$ and $k\in \mathbb N$; for simplicity of notation, we omit these indices. Moreover, the  Monge-Amp\`ere equation has the form  
	 \begin{equation}\label{Monge-Ampere-1}
	 	f(x)^{p^\star}= \frac{G_{\lambda,k}(T(x))}{\mathcal I_{\lambda,k}}{\rm det}DT(x)\ {\rm for \ a.e.}\ x\in \Omega.
	 \end{equation}
	Raising the latter equation to the power $1-\frac{1}{n}$, a reorganization and integration of the terms give 
	$$\int_\Omega\left(\frac{G_{\lambda,k}(T(x))}{\mathcal I_{\lambda,k}}\right)^{1-\frac{1}{n}}{\rm det}DT(x){\rm d}v_g(x)=\int_\Omega f(x)^{p^\star(1-\frac{1}{n})}({\rm det}DT(x))^\frac{1}{n}{\rm d}v_g(x).$$
	By a change of variables on the left hand side, see \cite[Corollary 4.7]{CEMS}, and by using Propositions \ref{Wang-Zhang-proposition} and \ref{Laplace-singular} for the right hand side, it follows that 
	\begin{eqnarray}\label{intermediate-ineq}
		\nonumber	\mathcal I_{\lambda,k}^{\frac{1}{n}-1}\int_M G_{\lambda,k}(y)^{1-\frac{1}{n}}{\rm d}v_g(y)&\leq& \int_\Omega f(x)^{p^\star(1-\frac{1}{n})}\left(1-\frac{\Delta_g u(x)}{n}\right){\rm d}v_g(x)\\&\leq& \int_\Omega f(x)^{p^\star(1-\frac{1}{n})}\left(1-\frac{\Delta_{g,\mathcal D'} u(x)}{n}\right){\rm d}v_g(x).
	\end{eqnarray}
	We now focus on the last term of \eqref{intermediate-ineq}. By  the divergence theorem and Schwarz inequality we have 
		\begin{eqnarray*}
			RHS&:= &\int_\Omega f(x)^{p^\star(1-\frac{1}{n})}\left(1-\frac{\Delta_{g,\mathcal D'} u(x)}{n}\right){\rm d}v_g(x)\\&=&\int_\Omega f^{p^\star(1-\frac{1}{n})}{\rm d}v_g+\frac{p^\star}{n}\left(1-\frac{1}{n}\right)\int_\Omega f(x)^{p^\star(1-\frac{1}{n})-1}{\nabla_g u(x)}\cdot \nabla_g f(x){\rm d}v_g(x)\\&\leq &\int_\Omega f^{p^\star(1-\frac{1}{n})}{\rm d}v_g+\frac{p^\star}{n}\left(1-\frac{1}{n}\right)\int_\Omega f(x)^{p^\star(1-\frac{1}{n})-1}|{\nabla_g u(x)}| \ | \nabla_g f(x)|{\rm d}v_g(x).
		\end{eqnarray*}
Here, $'\cdot'$ means the Riemannian inner product in $T_xM,$ $x\in M$.  Note that by \eqref{Optim-map} one has  for a.e.\ $x\in \Omega$  that $$|{\nabla_g u(x)}|	=d_g(x,T(x))\leq d_g(x,x_0)+d_g(x_0,T(x))\leq K_0+d_g(x_0,T(x)),$$	
where $K_0:=\max_{x\in \overline\Omega}d_g(x,x_0)<\infty$. Therefore, if we combine the latter two estimates, it follows that 
$$RHS\leq C(f,p,n,\Omega)+ \frac{p^\star}{n}\left(1-\frac{1}{n}\right)\int_\Omega f(x)^{p^\star(1-\frac{1}{n})-1}d_g(x_0,T(x)) | \nabla_g f(x)|{\rm d}v_g(x),$$
where 
$$C(f,p,n,\Omega)=\int_\Omega f(x)^{p^\star(1-\frac{1}{n})}{\rm d}v_g(x)+ K_0 \frac{p^\star}{n}\left(1-\frac{1}{n}\right)\int_\Omega f(x)^{p^\star(1-\frac{1}{n})-1} \ | \nabla_g f(x)|{\rm d}v_g(x)<\infty.$$
We emphasize that $C(f,p,n,\Omega)$ does \textit{not} depend either on  $\lambda>0$ or $k\in \mathbb N$. 

By applying H\"older's inequality and taking into account that  $\left(p^\star(1-\frac{1}{n})-1\right)p'=p^\star$, we obtain 
	\begin{eqnarray*}
	RHS&\leq &C(f,p,n,\Omega)+ \frac{p^\star}{n}\left(1-\frac{1}{n}\right)\left(\int_\Omega f(x)^{p^\star}d_g^{p'}(x_0,T(x)) {\rm d}v_g(x)\right)^\frac{1}{p'}\left(\int_\Omega  | \nabla_g f(x)|^p{\rm d}v_g(x)\right)^\frac{1}{p}.
	\end{eqnarray*}		
	The Monge-Amp\`ere equation \eqref{Monge-Ampere-1} and a change of variables imply  that
	\begin{eqnarray*}
	\int_\Omega f(x)^{p^\star}d_g^{p'}(x_0,T(x)) {\rm d}v_g(x)&=&\int_\Omega \frac{G_{\lambda,k}(T(x))}{\mathcal I_{\lambda,k}}d_g^{p'}(x_0,T(x)) {\rm det}DT(x){\rm d}v_g(x)\\&=&\int_M \frac{G_{\lambda,k}(y)}{\mathcal I_{\lambda,k}}d_g^{p'}(x_0,y) {\rm d}v_g(y);
	\end{eqnarray*}
	in particular, the latter term becomes independent on the optimal transport map $T$. 
	
	Summing up the above estimates, instead of the inequality \eqref{intermediate-ineq} we obtain 
\begin{eqnarray}\label{summ-ineq-0}
	\nonumber	\mathcal \mathcal I_{\lambda,k}^{\frac{1}{n}-\frac{1}{p}}\int_M G_{\lambda,k}(y)^{1-\frac{1}{n}}{\rm d}v_g(y)&\leq &C(f,p,n,\Omega)\mathcal I_{\lambda,k}^\frac{1}{p'}\\&&+\nonumber \frac{p^\star(n-1)}{n^2}\left(\int_M {G_{\lambda,k}(y)}d_g^{p'}(x_0,y) {\rm d}v_g(y)\right)^\frac{1}{p'}\left(\int_\Omega  | \nabla_g f(x)|^p{\rm d}v_g(x)\right)^\frac{1}{p}.
\end{eqnarray}

For every $\lambda>0$, we consider the Talentian bubble 
	$$
		G_{\lambda}(x)=\left(\lambda+d_g^{p'}(x_0,x)\right)^{-n},\ \ x\in M. 
	$$
	Note that $G_{\lambda,k}$ from \eqref{G-extrem-approx} has the property  $G_{\lambda,k}\leq G_{\lambda,k+1}\leq G_{\lambda}$ for every $\lambda>0$ and $k\in \mathbb N,$ and   $\lim_{k\to \infty}G_{\lambda,k}(x)=G_{\lambda}(x)$
	for every $x\in M$. Letting $k\to \infty $ in the  latter estimate, the monotone convergence theorem  implies that 
	\begin{eqnarray}\label{summ-ineq}
		\nonumber	\mathcal I_{\lambda}^{\frac{1}{n}-\frac{1}{p}}\int_M G_{\lambda}^{1-\frac{1}{n}}{\rm d}v_g&\leq &C(f,p,n,\Omega)\mathcal I_{\lambda}^\frac{1}{p'}\\&&+ \frac{p^\star(n-1)}{n^2}\left(\int_M {G_{\lambda}(y)}d_g^{p'}(x_0,y) {\rm d}v_g(y)\right)^\frac{1}{p'}\left(\int_M  | \nabla_g f|^p{\rm d}v_g\right)^\frac{1}{p},
	\end{eqnarray}
	where 
	$$\mathcal I_{\lambda}=\int_{M}G_{\lambda}(y){\rm d}v_g(y).$$
	
	We observe that every integral term (except the gradient) in \eqref{summ-ineq} can be written by means of the  function $\mathcal H$ defined in \eqref{H-function-definition}.  Indeed, we have
	$$\mathcal I_{\lambda}=\int_{M}G_{\lambda}{\rm d}v_g=\mathcal H(\lambda,n),\ \ \int_M G_{\lambda}^{1-\frac{1}{n}}{\rm d}v_g=\mathcal H(\lambda,n-1)$$
	and 
	$$\int_M {G_{\lambda}(y)}d_g^{p'}(x_0,y) {\rm d}v_g(y)=\mathcal H(\lambda,n-1)-\lambda \mathcal H(\lambda,n).$$
	Since $n>n-1>\frac{n}{p'}$, Proposition \ref{AT-estimates-proposition}/(i) implies that  these integrals are finite. Therefore, we have 
		\begin{eqnarray}\label{summ-ineq-2}
		\nonumber	\mathcal H(\lambda,n)^{\frac{1}{n}-\frac{1}{p}}\mathcal H(\lambda,n-1)&\leq &C(f,p,n,\Omega) \mathcal H(\lambda,n)^\frac{1}{p'}\\&&+ \frac{p^\star(n-1)}{n^2}\left(\mathcal H(\lambda,n-1)-\lambda \mathcal H(\lambda,n)\right)^\frac{1}{p'}\left(\int_M  | \nabla_g f|^p{\rm d}v_g\right)^\frac{1}{p}.
	\end{eqnarray}
	Let us multiply \eqref{summ-ineq-2} by $\lambda^{(n-1-\frac{n}{p'})\frac{1}{p'}}=\lambda^{\frac{n}{p^\star p'}}>0$; letting $\lambda\to \infty$,  the asymptotic properties from  Proposition \ref{AT-estimates-proposition}/(ii) (applied to $s=n$ and $s=n-1$, respectively) and an elementary computation imply that 	$$1 \leq
		{\sf AT}(n,p)\, {\sf AVR}_g^{-\frac{1}{n}}
		\left(\int_{M} |\nabla_g f|^p {\rm d}v_g\right)^\frac{1}{p},
$$
which is precisely inequality \eqref{egyenlet-0} in the normalized case \eqref{normalized-1}. In the generic case, we consider the function $\tilde f=\frac{\displaystyle|f|}{\left(\displaystyle\int_{M} |f|^{p^\star}{\rm d}v_g\right)^{{1}/{p^\star}}}$. 

It remains to prove the sharpness of the constant ${\sf AT}(n,p)\, {\sf AVR}_g^{-\frac{1}{n}}$ in \eqref{egyenlet-0}. 
Let us assume that there exists a constant $C>0$ such that $C<{\sf AT}(n,p)\, {\sf AVR}_g^{-\frac{1}{n}}$ and 
 for every $ f\in C_0^\infty(M)$ one has
\begin{equation}\label{egyenlet-2}
	\displaystyle  { \left(\int_{M} |f|^{p^\star}{\rm d}v_g\right)^{1/p^\star}\leq
		C
		\left(\int_{M} |\nabla_g f|^p {\rm d}v_g\right)^{1/p}}.
\end{equation}
Let $x_0\in M$ be fixed. For every $\lambda>0$, we can use the Talentian bubbles $M\ni x\mapsto (\lambda+d_g(x_0,x)^{p'})^\frac{p-n}{p}$  as test functions in   \eqref{egyenlet-2}; indeed, these functions may be approximated by elements from  $C_0^\infty(M)$ and standard limiting procedure yields that for every $\lambda>0$ one has
\begin{equation}\label{egyenlet-4}
\displaystyle  { \left(\int_{M} (\lambda+d_g(x_0,x)^{p'})^{-n}{\rm d}v_g\right)^{1/p^\star}\leq
	C\frac{n-p}{p}p'
	\left(\int_{M} (\lambda+d_g(x_0,x)^{p'})^{-n}d_g(x_0,x)^{p'} {\rm d}v_g\right)^{1/p}}.
\end{equation}
Here, we used the eikonal equation $|\nabla_g d_g(x_0,x)|=1$ for every $x\notin {\sf cut}(x_0).$
We notice that both integrals in \eqref{egyenlet-4} are finite, and they can be written by means of the function $\mathcal H$ from  \eqref{H-function-definition}; namely, it follows that
 \begin{equation}\label{egyenlet-5}
 	\displaystyle  {  \mathcal H(\lambda,n)^{\frac{1}{p^\star}}\leq
 		C\frac{n-p}{p}p'
 		\left(\mathcal H(\lambda,n-1)-\lambda \mathcal H(\lambda,n) \right)^\frac{1}{p}}.
 \end{equation}
If we multiply \eqref{egyenlet-5} by $\lambda^\frac{n}{p^\star p}>0$ and we take the limit $\lambda\to \infty$,  the asymptotic properties from  Proposition \ref{AT-estimates-proposition}/(ii) and a simple computation imply that  ${\sf AT}(n,p)\, {\sf AVR}_g^{-\frac{1}{n}}\leq C$, which contradicts our initial assumption $C<{\sf AT}(n,p)\, {\sf AVR}_g^{-\frac{1}{n}}$. The proof is concluded. \hfill $\square$

\begin{remark}\rm One of the crucial steps in the proof is the estimate after inequality  \eqref{intermediate-ineq}.  
Note that in the Euclidean setting, if we consider $\phi(x)=\frac{|x|^2}{2}-u(x)$, then one can write $$1-\frac{\Delta u(x)}{n}=\frac{\Delta \phi(x)}{n},$$ thus the classical divergence theorem can be applied for the \textit{whole} right hand side of \eqref{intermediate-ineq}, see \cite{CE-N-Villani}. However, in the Riemannian setting, we have to \textit{split} the terms in the right hand side of \eqref{intermediate-ineq}, thus beside of the expression containing the Laplacian -- where we  apply the divergence theorem -- we have to handle the extra term  
\begin{equation}\label{f-control}
	\int_\Omega f(x)^{p^\star(1-\frac{1}{n})}{\rm d}v_g(x),
\end{equation}
whose effect should \textit{not} be significant in our further limiting process. A reasonable way to guarantee the latter fact is to work with functions $f\in C_0^\infty(M)$. Although it is in principle possible to extend the class  $C_0^\infty(M)$ to appropriate Sobolev-regular functions,  possibly with non-compact support, and using OT theory for measures that do not necessarily have compact support, see Fathi and Figalli \cite{FF}, Figalli \cite{Figalli} or McCann \cite{McCann-Duke} (e.g., by requiring the finiteness of $p$-moments of the densities), we face serious technical difficulties. Indeed, beside a suitable control of the term  \eqref{f-control}, the applicability of a divergence-type theorem for non-smooth functions with not necessarily compact support poses another technical challenge in the present geometric setting, where no linear structure is available, which was crucial in the Euclidean framework,  see \cite[Lemma 7]{CE-N-Villani}. 
\end{remark}

\begin{remark}\label{remark-ledoux} \rm 
	The proof of the sharpness of the Sobolev constant in  \eqref{egyenlet-0} 
	provides, inter alia,  a new, elementary proof of the main result of Ledoux \cite{Ledoux}. In fact, this new approach -- based on the subtle asymptotic property of the Talentian bubbles generating the constant ${\sf AVR}_g$ -- will be applied in \S\ref{section-5} to prove the more general non-collapsing result of do Carmo and Xia \cite{doCarmo-Xia}. 
\end{remark}

%\vspace{5cm}
%
%On account of  \cite[Theorem 1.1]{Kristaly-Calculus} (see also  do Carmo and Xia \cite{doCarmo-Xia} for $p=2$), the valability of \eqref{egyenlet-2}
% implies that 
%$${\rm Vol }_g(B_x(r))	\geq \left(\frac{{\sf AT}(n,p)}{C}\right)^n\omega_n r^n,\ \ \forall r>0.$$
%In particular, by the latter estimate one has
%$$	{\sf AVR}_g\geq \left(\frac{{\sf AT}(n,p)}{C}\right)^n,$$
%which contradicts our assumption $C<{\sf AT}(n,p)\, {\sf AVR}_g^{-\frac{1}{n}}$. 

\subsection{The case $p=1$}
Let 	$ f\in C_0^\infty(M)$; as before, we may assume  that $f$ is nonnegative and 
\begin{equation}\label{normalized-2} \int_{M} f^{\frac{n}{n-1}}{\rm d}v_g=1.
\end{equation}
Let $\Omega=\{x\in M:f(x)>0\}$;  then $\overline \Omega$ is compact. 

Let $x_0\in \Omega$ and fix $\lambda>0$. We consider   the probability measures ${\rm d}\mu(x)= f(x)^{\frac{n}{n-1}} {\rm d}v_g(x)$ and ${\rm d}\nu(y)= \frac{\mathbbm{1}_{B_{x_0}(\lambda)}(y)}{{\rm Vol}_g(B_{x_0}(\lambda))}{\rm d}
v_g(y)$  on $(M,g)$, where $\mathbbm{1}_S$ stands for the characteristic function of the set $S\subset M.$  According to McCann \cite{McCann}, there exists a unique  optimal transport
map $T:\overline \Omega\to \overline {B_{x_0}(\lambda)}\subset M$  pushing $\mu$ forward to $\nu$ with the form 
$
	T(x)=\exp_x(-\nabla_g u(x))$  for a.e.\  $x\in \overline \Omega,
$
where $u:\overline \Omega\to \mathbb R$  is a $c=d_g^2/2$-concave function. Moreover, the  Monge-Amp\`ere equation reads as
\begin{equation}\label{Monge-Ampere-2}
	f(x)^{\frac{n}{n-1}}= \frac{1}{{\rm Vol}_g(B_{x_0}(\lambda))}{\rm det}DT(x)\ {\rm for \ a.e.}\ x\in \Omega.
\end{equation}
A simple rearrangement of \eqref{Monge-Ampere-2} implies
$$\int_\Omega\left(\frac{1}{{\rm Vol}_g(B_{x_0}(\lambda))}\right)^{1-\frac{1}{n}}{\rm det}DT(x){\rm d}v_g(x)=\int_\Omega f(x)({\rm det}DT(x))^\frac{1}{n}{\rm d}v_g(x).$$
A change of variables at the left hand side, combined with Propositions \ref{Wang-Zhang-proposition} and \ref{Laplace-singular}, and the divergence theorem for the right hand side  yield 
\begin{eqnarray}\label{isso-estimate}
\nonumber	{\rm Vol}_g(B_{x_0}(\lambda))^\frac{1}{n}&\leq& \int_\Omega f(x)\left(1-\frac{\Delta_{g,\mathcal D'} u(x)}{n}\right){\rm d}v_g(x)=\int_\Omega f{\rm d}v_g+\frac{1}{n}\int_\Omega \nabla_g f(x)\cdot \nabla_g u(x){\rm d}v_g(x)\\&\leq &\int_\Omega f{\rm d}v_g+\frac{1}{n}\|\nabla_g u\|_{{L^\infty}(\Omega)}\int_\Omega |\nabla_g f| {\rm d}v_g.
\end{eqnarray}
Since $T(x)\in \overline {B_{x_0}(\lambda)}$ for a.e.\ $x\in \Omega$,  we have $$|{\nabla_g u(x)}|	=d_g(x,T(x))\leq d_g(x,x_0)+d_g(x_0,T(x))\leq K_0+\lambda,$$	
where $K_0:=\max_{x\in \overline\Omega}d_g(x,x_0)<\infty$. This estimate and \eqref{isso-estimate} imply that 
$${\rm Vol}_g(B_{x_0}(\lambda))^\frac{1}{n}\leq \int_\Omega f{\rm d}v_g+\frac{1}{n}(K_0+\lambda)\int_\Omega |\nabla_g f| {\rm d}v_g.$$
Dividing by $\lambda>0$ the latter inequality and letting $\lambda\to \infty$, it follows that 
$$(\omega_n{\sf AVR}_g)^{\frac{1}{n}}=\lim_{\lambda\to \infty}\frac{{\rm Vol}_g(B_{x_0}(\lambda))^\frac{1}{n}}{\lambda}\leq \frac{1}{n}\int_\Omega |\nabla_g f| {\rm d}v_g,$$
which is precisely \eqref{egyenlet-0} for $p=1$ due to  \eqref{normalized-2}; the non-normalized case follows in a standard way. 

Note that the Sobolev inequality \eqref{egyenlet-0} for $p=1$ is equivalent to the  isoperimetric inequality 
\begin{equation}\label{isoperimetrikus}
	n(\omega_n{\sf AVR}_g)^{\frac{1}{n}}{\rm Vol}_g(\Omega)^\frac{n-1}{n}\leq \mathcal P_g(\partial \Omega)
\end{equation}
for every bounded open domain $\Omega\subset M$ with smooth boundary ($\mathcal P_g$ being the perimeter); this equivalence is rather standard and can be found for instance in Ambrosio, Carlotto and Massaccesi \cite[Lemma A.24]{ACM}.  On the other hand, we know that  \eqref{isoperimetrikus} is sharp, see Balogh and Krist\'aly \cite{BK} and Brendle \cite{Brendle}; thus,  it turns out that  \eqref{egyenlet-0} is also sharp.  \hfill $\square$

\section{Proof of the sharp $L^p$-logarithmic Sobolev inequality (Theorem \ref{log-Sobolev-main})}\label{section-4}
\subsection{The case $p>1$}

%\subsection{Proof of Theorem \ref{log-Sobolev-main}} 
Let $p>1$  and fix	$ f\in C_0^\infty(M)$ arbitrarily;  we may assume that $f$ is nonnegative and 
$$\int_{M} f^{p}{\rm d}v_g=1.$$
As before, let $\Omega=\{x\in M:f(x)>0\}$; since $ f\in C_0^\infty(M)$, then $\overline \Omega$ is compact. 

Let $x_0\in \Omega.$ For every $\lambda>0$ and $k\in \mathbb N$, we introduce the truncated Gaussian bubble  $G_{\lambda,k}:M\to \mathbb R$ given by
$$
	G_{\lambda,k}(x)=P_k(d_g(x_0,x))e^{-\lambda d_g^{p'}(x_0,x)},
$$
where $P_k$ is defined in \eqref{truncation-k}. We observe that the  support of $G_{\lambda,k}$ is the ball $\overline {B_{x_0}(k+1)}$. Let $$\mathcal J_{\lambda,k}=\int_{M}G_{\lambda,k}(y){\rm d}v_g(y);$$
clearly, $0<\mathcal J_{\lambda,k}<\infty$ for every $\lambda>0$ and $k\in \mathbb N$. 

Let  ${\rm d}\mu(x)= f^{p}(x) {\rm d}v_g(x)$ and ${\rm d}\nu(y)= \frac{G_{\lambda,k}(y)}{\mathcal J_{\lambda,k}}{\rm d}
v_g(y)$ be two probability measures on $(M,g)$ with compact supports; by the theory of OT one can find  a unique  
map $T:\overline \Omega\to \overline {B_{x_0}(k+1)}\subset M$  pushing $\mu$ forward to $\nu$ having the form 
$
	T(x)=\exp_x(-\nabla_g u(x))$  for  a.e.\ $  x\in \overline \Omega,
$
for some $c=d_g^2/2$-concave function $u:\overline \Omega\to \mathbb R$.  The associated Monge-Amp\`ere equation is
\begin{equation}\label{Monge-Ampere-log}
	f^{p}(x)= \frac{G_{\lambda,k}(T(x))}{\mathcal J_{\lambda,k}}{\rm det}DT(x)\ {\rm for \ a.e.}\ x\in \Omega.
\end{equation}
Accordingly, by \eqref{Monge-Ampere-log}, a change of variables, Jensen's inequality and Propositions \ref{Wang-Zhang-proposition} and \ref{Laplace-singular}, we have
\begin{eqnarray*}
	\int_{M}f^p\log f^p{\rm d}v_g &=&\int_{\Omega}f^p(x)\log \left(\frac{G_{\lambda,k}(T(x))}{\mathcal J_{\lambda,k}}{\rm det}DT(x)\right) {\rm d}v_g(x)\\&=&
	\int_{\Omega}f^p(x)\log \left({G_{\lambda,k}(T(x))}\right) {\rm d}v_g(x)+n\int_{\Omega}f^p(x)\log \left(\frac{{\rm det}^\frac{1}{n}DT(x)}{\mathcal J^\frac{1}{n}_{\lambda,k}}\right) {\rm d}v_g(x)\\&\leq&
\int_{\Omega}\frac{G_{\lambda,k}(T(x))}{\mathcal J_{\lambda,k}}\log \left({G_{\lambda,k}(T(x))}\right){\rm det}DT(x) {\rm d}v_g(x)\\&&+n\log\left(\int_{\Omega}f^p(x) \frac{{\rm det}^\frac{1}{n}DT(x)}{\mathcal J^\frac{1}{n}_{\lambda,k}}{\rm d}v_g(x)\right) \\&\leq&
\frac{1}{{\mathcal J_{\lambda,k}}}\int_{M}{G_{\lambda,k}(y)}\log \left({G_{\lambda,k}(y)}\right) {\rm d}v_g(y)\\&&+n\log\left(\int_{\Omega} \frac{f^p(x)}{\mathcal J^\frac{1}{n}_{\lambda,k}}\left(1-\frac{\Delta_g u(x)}{n}\right) {\rm d}v_g(x)\right)
\\&\leq&
\frac{1}{{\mathcal J_{\lambda,k}}}\int_{M}{G_{\lambda,k}}\log \left({G_{\lambda,k}}\right) {\rm d}v_g+n\log\left(\int_{\Omega} \frac{f^p(x)}{\mathcal J^\frac{1}{n}_{\lambda,k}}\left(1-\frac{\Delta_{g,\mathcal D'} u(x)}{n}\right) {\rm d}v_g(x)\right).
\end{eqnarray*}
On the other hand, by the divergence theorem and Schwarz inequality one has  that
\begin{eqnarray*}
	-\int_\Omega f^p(x)\Delta_{g,\mathcal D'} u(x){\rm d}v_g(x)&=&p\int_\Omega f^{p-1}(x)\nabla_g f(x)\cdot\nabla_g u(x){\rm d}v_g(x)\\&\leq& p\int_\Omega f^{p-1}(x)|\nabla_g f(x)|\ |\nabla_g u(x)|{\rm d}v_g(x)\\&\leq & K_0 p \int_\Omega f^{p-1} |\nabla_g f|{\rm d}v_g+ p\int_\Omega f^{p-1}(x)d_g(x_0,T(x)) |\nabla_g f(x)|{\rm d}v_g(x),
\end{eqnarray*}
where $K_0:=\max_{x\in \overline\Omega}d_g(x,x_0)<\infty$ and we applied for a.e.\ $x\in \Omega$ that  $$|{\nabla_g u(x)}|	=d_g(x,T(x))\leq d_g(x,x_0)+d_g(x_0,T(x))\leq K_0+d_g(x_0,T(x)).$$
By H\"older's inequality and   Monge-Amp\`ere equation \eqref{Monge-Ampere-log}, it turns out that 
\begin{eqnarray*}
	-\int_\Omega f^p(x)\Delta_{g,\mathcal D'} u(x){\rm d}v_g(x)&\leq & K_0 p \int_\Omega f^{p-1} |\nabla_g f|{\rm d}v_g\\&&+ p\left(\int_\Omega  |\nabla_g f|^p{\rm d}v_g\right)^\frac{1}{p}\left(\int_\Omega f^{p}(x)d_g^{p'}(x_0,T(x)) {\rm d}v_g(x)\right)^\frac{1}{p'}\\ &= & K_0 p \int_\Omega f^{p-1} |\nabla_g f|{\rm d}v_g\\&&+ p\left(\int_\Omega  |\nabla_g f|^p{\rm d}v_g\right)^\frac{1}{p}\left(\int_M \frac{G_{\lambda,k}(y)}{\mathcal J_{\lambda,k}}d_g^{p'}(x_0,y) {\rm d}v_g(y)\right)^\frac{1}{p'}. 
\end{eqnarray*}
The above estimates imply that
\begin{eqnarray*}
	\int_{M}f^p\log f^p{\rm d}v_g &\leq&\frac{1}{{\mathcal J_{\lambda,k}}}\int_{M}{G_{\lambda,k}}\log \left({G_{\lambda,k}}\right) {\rm d}v_g\\&&+n\log\left.\Biggl[\frac{1}{\mathcal J^\frac{1}{n}_{\lambda,k}} \left. \Biggl(C(f,p,\Omega)+\right.\right.\\&&\ \ \ \ \ \  \ \ \ \left.\left.+\frac{p}{n}\left(\int_\Omega  |\nabla_g f|^p{\rm d}v_g\right)^\frac{1}{p}\left(\int_M \frac{G_{\lambda,k}(y)}{\mathcal J_{\lambda,k}}d_g^{p'}(x_0,y) {\rm d}v_g(y)\right)^\frac{1}{p'}\right)\right.\Biggl],
\end{eqnarray*}
where 
$$C(f,p,\Omega)=1+ K_0 \frac{p}{n} \int_\Omega f^{p-1} |\nabla_g f|{\rm d}v_g.$$
On account of  the definition of functions $\mathcal L_i$, $i\in \{1,2\}$, see \eqref{LL-function-definition},  the monotone convergence theorem implies that for every $\lambda>0$ one has 
$$\lim_{k\to \infty}\mathcal J_{\lambda,k}= \mathcal L_1(\lambda)$$
and 
$$\lim_{k\to \infty}\int_{M}{G_{\lambda,k}}\log \left({G_{\lambda,k}}\right) {\rm d}v_g=-\lambda \mathcal L_2(\lambda)=-\lambda\lim_{k\to \infty}\int_{M}{G_{\lambda,k}(y)}d_g^{p'}(x_0,y) {\rm d}v_g(y).$$
Letting $k\to \infty$ in the latter inequality, we infer that
\begin{eqnarray}\label{4-4-ineq}
	\int_{M}f^p\log f^p{\rm d}v_g \nonumber &\leq&-\frac{\lambda\mathcal L_2(\lambda)}{\mathcal L_1(\lambda)}\\&&+n\log\left(\frac{1}{\mathcal L_1(\lambda)^\frac{1}{n}} \left(C(f,p,\Omega)+\frac{p}{n}\left(\int_\Omega  |\nabla_g f|^p{\rm d}v_g\right)^\frac{1}{p}\left(\frac{\mathcal L_2(\lambda)}{\mathcal L_1(\lambda)}\right)^\frac{1}{p'}\right)\right).
\end{eqnarray}
Due to Proposition \ref{L-estimates-proposition}/(ii), one has 
$$\lim_{\lambda\to 0} \frac{\lambda\mathcal L_2(\lambda)}{\mathcal L_1(\lambda)}=\frac{n}{p'},\ \ \lim_{\lambda\to 0} \frac{1}{\mathcal L_1(\lambda)}=0,\ \  {\rm and}\ \ \lim_{\lambda\to 0} \frac{\mathcal L_2(\lambda)^\frac{1}{p'}}{\mathcal L_1(\lambda)^{\frac{1}{p'}+\frac{1}{n}}}=\left(\frac{n}{p'}\right)^\frac{1}{p'}{\left(\omega_n{\sf AVR}_g\Gamma\left(\frac{n}{p'}+1\right)\right)^{-\frac{1}{n}}}.$$ 
Letting now $\lambda\to 0$ in \eqref{4-4-ineq}, the above limits provide 
$$
\int_{M}f^p\log f^p{\rm d}v_g\leq \frac{n}{p}\log\left({\sf L}({n,p}){\sf AVR}_ g^{-\frac{p}{n}}\int_M |\nabla_g f|^p{\rm d}v_g \right), 
$$
which is precisely inequality \eqref{LSI}. 

We now assume by contradiction that there exists a constant $C>0$ such that $C<{\sf L}(n,p)\, {\sf AVR}_g^{-\frac{p}{n}}$ and 
for every $ f\in C_0^\infty(M)$ with $\displaystyle\int_{M}|f|^p{\rm d}v_g=1$  one has
\begin{equation}\label{egyenlet-30}
	\int_{M}|f|^p\log |f|^p{\rm d}v_g\leq \frac{n}{p}\log\left(C\int_M |\nabla_g f|^p{\rm d}v_g \right). 
\end{equation}
If we fix $x_0\in M$, one can consider for every $\lambda>0$  the Gaussian bubbles $$M\ni x\mapsto {e^{-\frac{\lambda}{p} d_g^{p'}(x_0,x)}}{\mathcal L_1(\lambda)^{-\frac{1}{p}}}$$ as test functions in   \eqref{egyenlet-30} by using suitable approximation with  elements from  $C_0^\infty(M)$. Accordingly, inserting the above Gaussian bubble into \eqref{egyenlet-30}, we infer that for every $\lambda>0$ one has
$$\int_{M}\frac{e^{-{\lambda} d_g^{p'}(x_0,x)}}{\mathcal L_1(\lambda)}\log \left(\frac{e^{-{\lambda} d_g^{p'}(x_0,x)}}{\mathcal L_1(\lambda)}\right){\rm d}v_g\leq \frac{n}{p}\log\left(\frac{C}{\mathcal L_1(\lambda)}\left(\frac{\lambda p'}{p}\right)^p \int_Me^{-\lambda d_g^{p'}(x_0,x)}d_g^{p'}(x_0,x){\rm d}v_g(x) \right).$$
On account of \eqref{LL-function-definition}, the latter inequality is equivalent to 
$$-\frac{\lambda\mathcal L_2(\lambda)}{\mathcal L_1(\lambda)}\leq \frac{n}{p}\log\left(C\left(\frac{\lambda p'}{p}\right)^p \frac{\mathcal L_2(\lambda)}{\mathcal L_1(\lambda)^{1-\frac{p}{n}}} \right)
$$
for every $\lambda>0.$
By Proposition \ref{L-estimates-proposition}/(ii), if  $\lambda\to 0$ in the above inequality, it follows that
$$-\frac{n}{p'}\leq \frac{n}{p}\log\left(C\left(\frac{ p'}{p}\right)^p \frac{n}{p'}\left(\omega_n{\sf AVR}_g\Gamma\left(\frac{n}{p'}+1\right)\right)^\frac{p}{n} \right),$$
which is equivalent to ${\sf L}(n,p)\, {\sf AVR}_g^{-\frac{p}{n}}\leq C$; but this inequality contradicts our initial assumption $C<{\sf L}(n,p)\, {\sf AVR}_g^{-\frac{p}{n}}$. This ends the proof of the sharpness in  \eqref{LSI}. \hfill $\square$

%\vspace{5cm}
%Due to  \cite[Theorem 1.2]{Kristaly-Calculus}, by \eqref{egyenlet-3}
%one has the volume non-collapsing estimate
%$${\rm Vol }_g(B_{x_0}(r))	\geq \left(\frac{{\sf L}(n,p)}{C}\right)^\frac{n}{p}\omega_n r^n,\ \ \forall r>0;$$
%thus, one has that
%$$	{\sf AVR}_g\geq \left(\frac{{\sf L}(n,p)}{C}\right)^\frac{n}{p},$$
%contradicting the assumption $C<{\sf L}(n,p)\, {\sf AVR}_g^{-\frac{p}{n}}$. \hfill $\square$

\subsection{The case $p=1$}
Let 	$ f\in C_0^\infty(M)$ with $\displaystyle\int_{M} |f|{\rm d}v_g=1$. As in the previous cases, we may assume  that $f$ is nonnegative and let $\Omega=\{x\in M:f(x)>0\}$;  then $\overline \Omega$ is compact. 

Let $x_0\in \Omega$ and fix $\lambda>0$. Let  ${\rm d}\mu(x)= f(x) {\rm d}v_g(x)$ and ${\rm d}\nu(y)= \frac{\mathbbm{1}_{B_{x_0}(\lambda)}(y)}{{\rm Vol}_g(B_{x_0}(\lambda))}{\rm d}
v_g(y)$  be two probability measures; then there exists a unique  optimal transport
map $T:\overline \Omega\to \overline {B_{x_0}(\lambda)}\subset M$  pushing $\mu$ forward to $\nu$ having the form
$
	T(x)=\exp_x(-\nabla_g u(x))$   for  a.e.\ $  x\in \overline \Omega,
$
where $u:\overline \Omega\to \mathbb R$  is a $c=d_g^2/2$-concave function. In this case,  the  Monge-Amp\`ere equation becomes 
\begin{equation}\label{Monge-Ampere-23}
	f(x)= \frac{1}{{\rm Vol}_g(B_{x_0}(\lambda))}{\rm det}DT(x)\ {\rm for \ a.e.}\ x\in \Omega.
\end{equation}
As in the previous subsection,   by using relation \eqref{Monge-Ampere-23}, Jensen's inequality,  Propositions \ref{Wang-Zhang-proposition} and \ref{Laplace-singular} and the divergence theorem, we obtain that 
\begin{eqnarray*}
	\int_{M}f\log f{\rm d}v_g &=&\int_{\Omega}f(x)\log \left(\frac{1}{{\rm Vol}_g(B_{x_0}(\lambda))}{\rm det}DT(x)\right) {\rm d}v_g(x)\\&=&
	-\log({\rm Vol}_g(B_{x_0}(\lambda)))+n\int_{\Omega}f(x)\log \left({({\rm det}DT(x))^\frac{1}{n}}\right) {\rm d}v_g(x)\\&\leq&
		-\log({\rm Vol}_g(B_{x_0}(\lambda)))+n\log\left(\int_{\Omega}f(x) {({\rm det}DT(x))^\frac{1}{n}}{\rm d}v_g(x)\right) \\&\leq&
	-\log({\rm Vol}_g(B_{x_0}(\lambda)))+n\log\left(\int_{\Omega} {f(x)}\left(1-\frac{\Delta_{g,\mathcal D' }u(x)}{n}\right) {\rm d}v_g(x)\right)\\&=&
	-\log({\rm Vol}_g(B_{x_0}(\lambda)))+n\log\left(1+\frac{1}{n}\int_{\Omega} {\nabla_gf(x)}\cdot \nabla_g u(x){\rm d}v_g(x)\right) .
\end{eqnarray*}
If $K_0:=\max_{x\in \overline\Omega}d_g(x,x_0)<\infty$, then for a.e.\ $x\in \Omega$ one has  $|{\nabla_g u(x)}|	=d_g(x,T(x))\leq d_g(x,x_0)+d_g(x_0,T(x))\leq  K_0+\lambda;$
 therefore,  
$$\int_{\Omega} {\nabla_gf(x)}\cdot \nabla_g u(x){\rm d}v_g(x)\leq \int_{\Omega} |{\nabla_gf(x)}|\ | \nabla_g u(x)|{\rm d}v_g(x)\leq (K_0+\lambda) \int_{\Omega} |{\nabla_gf}|{\rm d}v_g.$$
Combining the above estimates, it follows for every $\lambda>0$ that
\begin{eqnarray*}
	\int_{M}f\log f{\rm d}v_g&\leq &-\log({\rm Vol}_g(B_{x_0}(\lambda)))+n\log\left(1+\frac{K_0+\lambda}{n} \int_{\Omega} |{\nabla_gf}|{\rm d}v_g\right)\\&=&n\log\left(\frac{1}{{\rm Vol}^\frac{1}{n}_g(B_{x_0}(\lambda))}+\frac{K_0+\lambda}{n{\rm Vol}^\frac{1}{n}_g(B_{x_0}(\lambda))} \int_{\Omega} |{\nabla_gf}|{\rm d}v_g\right).
\end{eqnarray*}
Letting $\lambda\to \infty$, we infer that
$$\int_{M}f\log f{\rm d}v_g\leq n\log\left(\frac{1}{n(\omega_n {\sf AVR}_g)^\frac{1}{n}} \int_{\Omega} |{\nabla_gf}|{\rm d}v_g\right),$$
which is precisely inequality \eqref{LSI} for $p=1.$

Concerning the sharpness of the constant in \eqref{LSI} for $p=1$, we assume there exists $C<n^{-1}(\omega_n {\sf AVR}_g)^{-\frac{1}{n}}={\sf L}({n,1}){\sf AVR}_g^{-\frac{1}{n}}$ such that for every $ f\in C_0^\infty(M)$ with $\displaystyle\int_{M} |f|{\rm d}v_g=1$, 
$$\int_{M}f\log f{\rm d}v_g\leq n\log\left(C \int_{\Omega} |{\nabla_gf}|{\rm d}v_g\right).$$
By density argument, the above inequality is valid also for functions belonging to ${\rm BV}(M)$ (with the usual modification of the right hand side); in particular, we may consider $f:=\frac{\mathbbm{1}_\Omega}{{\rm Vol}_g(\Omega)}$ for every bounded open set $\Omega\subset M$ with smooth boundary. With this choice the latter inequality reduces to $$-\log({\rm Vol}_g(\Omega))\leq n\log\left(C \frac{\mathcal P_g(\partial\Omega)}{{\rm Vol}_g(\Omega)}\right),$$
which is equivalent to the isoperimetric inequality
$$C^{-1}{\rm Vol}_g(\Omega)^\frac{n-1}{n}\leq \mathcal P_g(\partial\Omega),$$
valid for every bounded open set $\Omega\subset M$ with smooth boundary. However, the initial assumption $n(\omega_n {\sf AVR}_g)^{\frac{1}{n}}<C^{-1}$ contradicts the sharpness of the isoperimetric inequality \eqref{isoperimetrikus}. \hfill $\square$

\subsection{Application: Gaussian logarithmic-Sobolev inequality (Theorem \ref{Gaussian-log-Sobolev-main})} Let $ h\in C_0^\infty(M)$  be such that
 $\displaystyle\int_M h^2 {\rm d}\gamma_V=1$. We consider the function $f=\sqrt{\rho} h,$
 with $\rho=G_V^{-1}e^{-V}$. It is clear that
 $$\int_M f^2{\rm d}v_g=\int_M  h^2\rho{\rm d}v_g=\int_M  h^2{\rm d}\gamma_V=1.$$
 In particular, since ${\sf L}({n,2})=\frac{2}{\pi n e},$
 Theorem \ref{log-Sobolev-main} implies (for $p=2$) that
  	\begin{equation}\label{LSI=p=2}
  	\int_{M}f^2\log f^2{\rm d}v_g\leq \frac{n}{2}\log\left(\frac{2}{\pi n e}{\sf AVR}_ g^{-\frac{2}{n}}\int_M |\nabla_g f|^2{\rm d}v_g \right). 
  \end{equation}
First, we have that
\begin{equation}\label{entropy-term}
	\int_{M}f^2\log f^2{\rm d}v_g=\int_{M}h^2\log h^2{\rm d}\gamma_V-\log G_V-\int_M h^2V{\rm d}\gamma_V.
\end{equation}
Moreover, a simple computation also yields that
\begin{eqnarray*}
	|\nabla_g f|^2=\rho\left|\nabla_gh-\frac{1}{2}h\nabla_g V\right|^2=\rho\left|\nabla_gh\right|^2+\frac{G_V^{-1}}{2}\nabla_g h^2 \cdot \nabla_g e^{-V}+\frac{\rho}{4}h^2|\nabla_gV|^2.
\end{eqnarray*}
Therefore, by the divergence theorem one has
\begin{eqnarray*}
\int_M	|\nabla_g f|^2{\rm d}v_g&=&\int_M\left|\nabla_gh\right|^2{\rm d}\gamma_V+\frac{G_V^{-1}}{2}\int_M\nabla_g h^2 \cdot \nabla_g e^{-V}{\rm d}v_g+\frac{G_V^{-1}}{4}\int_Mh^2|\nabla_gV|^2e^{-V}{\rm d}v_g\\&=&\int_M\left|\nabla_gh\right|^2{\rm d}\gamma_V+\frac{G_V^{-1}}{2}\int_M h^2 \left(-\Delta_g (e^{-V})+\frac{1}{2}|\nabla_gV|^2e^{-V}\right){\rm d}v_g\\&=&\int_M\left|\nabla_gh\right|^2{\rm d}\gamma_V+\frac{G_V^{-1}}{2}\int_M h^2 \left(-\frac{1}{2}|\nabla_gV|^2+\Delta_g V\right)e^{-V}{\rm d}v_g.
\end{eqnarray*}
Combining the inequality $\log (et)\leq t$ for every $t>0$ with the latter computation, it turns out that 
\begin{eqnarray*}
	\log\left(\frac{2}{\pi n e}{\sf AVR}_ g^{-\frac{2}{n}}\int_M |\nabla_g f|^2{\rm d}v_g \right)&=&\log\left(\frac{1}{2\pi  e^2}{\sf AVR}_ g^{-\frac{2}{n}} \right)+\log\left(\frac{4e}{ n }\int_M |\nabla_g f|^2{\rm d}v_g \right)\\&\leq &
	\log\left(\frac{1}{2\pi  e^2}{\sf AVR}_ g^{-\frac{2}{n}} \right)+\frac{4}{ n }\int_M |\nabla_g f|^2{\rm d}v_g \\&\leq &
	\log\left(\frac{1}{2\pi  e^2}{\sf AVR}_ g^{-\frac{2}{n}} \right)+\frac{4}{ n }\int_M |\nabla_g h|^2{\rm d}\gamma_V\\&&+ \frac{2}{n}\int_M h^2 \left(-\frac{1}{2}|\nabla_gV|^2+\Delta_g V\right){\rm d}\gamma_V.
\end{eqnarray*}
Replacing \eqref{entropy-term} and the latter estimate into  
\eqref{LSI=p=2}, it follows by our assumption \eqref{assumption-Cv} that 
\begin{eqnarray*}
	\int_{M}h^2\log h^2{\rm d}\gamma_V&\leq & 2\int_M |\nabla_g h|^2{\rm d}\gamma_V+\log\left(\frac{G_V}{(2\pi)^\frac{n}{2} {\sf AVR}_g}\right)\\&&+\int_M h^2 \left(V-\frac{1}{2}|\nabla_gV|^2+\Delta_g V-n\right){\rm d}\gamma_V\\&\leq&
	2\int_M |\nabla_g h|^2{\rm d}\gamma_V+\log\left(\frac{G_Ve^{C_V}}{(2\pi)^\frac{n}{2} {\sf AVR}_g}\right),
\end{eqnarray*}
which is precisely the desired inequality \eqref{Gaussian-LSI}. \hfill $\square$ 

\begin{remark}\label{remark-particular}\rm 
 (a) If $V(x)=\frac{1}{2}d_g^2(x_0,x)$ for some $x_0\in M$, then by the eikonal equation $|\nabla_g d_g(x_0,\cdot)|=1$ ${\rm d}v_g$-a.e. on $M$ and the Laplace comparison $\Delta_g V\leq n$ (note that ${\sf Ric}\geq 0$), it turns out that assumption  \eqref{assumption-Cv} holds with the choice $C_V=0.$ Moreover, by a similar computation as in 
Proposition \ref{L-estimates-proposition}, combined with the Bishop-Gromov volume comparison principle gives	$$G_V=\int_M e^{-V}{\rm d}v_g \leq (2\pi)^\frac{n}{2}.$$
Therefore, by inequality \eqref{Gaussian-LSI} we obtain 
$$
	\int_{M}h^2\log h^2{\rm d}\gamma_V\leq 2\int_M |\nabla_g h|^2{\rm d}\gamma_V-\log{\sf AVR}_g,
$$
which is exactly \eqref{Gaussian-LSI-particular}. 

(b) In general, if $V(x)=s(d_g(x_0,x))$ for some nonnegative and non-decreasing function $s\in C^2(0,\infty)$, then we may replace assumption \eqref{assumption-Cv} by
$$s(t)-\frac{s'(t)^2}{2}+s'(t)\frac{n-1}{t}+s''(t)-n\leq C_V\ \ {\rm a.e.}\ t>0.$$

(c) A closer look to the proof of Theorem \ref{Gaussian-log-Sobolev-main} provides a slightly more general form of the Gaussian logarithmic-Sobolev inequality. Indeed, if we assume instead of \eqref{assumption-Cv} that there exists $K>0$ and $C_V\in \mathbb R$ such that 
\begin{equation}\label{assumption-Cv-modified}
	V-\frac{|\nabla_g V|^2}{2K}+\frac{\Delta_gV}{K} -n\leq C_V\ \ {\rm a.e.\ on} \ M,
\end{equation}
then for every $ h\in C_0^\infty(M)$ 
with $\displaystyle\int_M h^2 {\rm d}\gamma_V=1$, one has that 
$$
	\int_{M}h^2\log h^2{\rm d}\gamma_V\leq \frac{2}{K}\int_M |\nabla_g h|^2{\rm d}\gamma_V+\log\left(\frac{K^\frac{n}{2} G_Ve^{C_V}}{(2\pi)^\frac{n}{2} {\sf AVR}_g}\right).
$$
Note that $V_K(x)=\frac{K}{2}d_g^2(x_0,x)$ satisfies \eqref{assumption-Cv-modified} with $C_{V_K}=0$; moreover, $G_{V_K}=\displaystyle\int_M e^{-V_K}{\rm d}v_g \leq (2\pi/K)^\frac{n}{2}.$ Thus, in this particular case, we have for every $ h\in C_0^\infty(M)$ 
with $\displaystyle\int_M h^2 {\rm d}\gamma_{V_K}=1$  that
$$	\int_{M}h^2\log h^2{\rm d}\gamma_{V_K}\leq \frac{2}{K}\int_M |\nabla_g h|^2{\rm d}\gamma_{V_K}-\log{\sf AVR}_g.$$
This inequality is well-known in the Euclidean setting (${\sf AVR}_g=1$), which is the dimension-free logarithmic-Sobolev inequality of Gross \cite{Gross}. 
\end{remark}

\section{Byproduct: a simple proof of do Carmo and Xia's non-collapsing result}\label{section-5}

The purpose of the present section is to give a short proof of the main result of do Carmo and Xia \cite{doCarmo-Xia}, which deals with a weighted Sobolev (or, Caffarelli--Kohn--Nirenberg) inequality on Riemannian manifolds with ${\sf Ric}\geq 0.$

Let $n\geq 3$ be an integer,  the numbers $a,b\in \mathbb R$ such that 
\begin{equation}\label{a-b}
	0\leq a<\frac{n-2}{2},\ \ a\leq b<a+1,
\end{equation}
and set the constants 
$$q=\frac{2n}{n-2+2(b-a)}\ \ {\rm and} \ \  K_{a,b}=((n-2a-2)(n-bq))^{-\frac{1}{2}}\left(\frac{2-bq+2a}{n\omega_n}\frac{\Gamma\left(\frac{n}{a+1-b}\right)}{\Gamma^2\left(\frac{n}{2(a+1-b)}\right)}\right)^\frac{a+1-b}{n}.$$
Chou and  Chu \cite{Chou} proved that for every $f\in C_0^\infty(\mathbb R^n)$ one has
$$\left(\int_{\mathbb R^n}|x|^{-bq}|f|^q{\rm d}x\right)^\frac{1}{q}\leq K_{a,b} \left(\int_{\mathbb R^n}|x|^{-2a}|\nabla f|^2{\rm d}x\right)^\frac{1}{2},$$
and  $K_{a,b}$ is optimal, the class of extremal functions $(f_\lambda)_{\lambda>0}$  are Talenti bubbles of the form 
$$f_\lambda(x)=\left(\lambda+|x|^{2-bq+2a}\right)^{-\frac{n-2a-2}{2-bq+2a}},\ \ x\in \mathbb R^n.$$
Note that if $a=b=0$, then $K_{0,0}={\sf AT}(n,2)$ and the latter inequality reduces to the sharp $L^2$-Sobolev inequality \eqref{Sobolev-000}. 

The main non-collapsing rigidity result of do Carmo and Xia \cite[Theorem 1.1]{doCarmo-Xia} reads as follows: 

\begin{theorem}\label{thm-docarmo}
	Let $C\geq K_{a,b}$ be a constant, $(M,g)$ be a noncompact, complete  $n$-dimensional Riemannian manifold $(n\geq 3)$ with ${\sf Ric}\geq 0$, and the constants $a,b\in \mathbb R$ verifying relations \eqref{a-b}. If $x_0\in M$ is fixed and for every $f\in C_0^\infty(M)$ one has
	\begin{equation}\label{egyenlet-ujra}
		\displaystyle  { \left(\int_{M} d_g^{-bq}(x_0,x)|f(x)|^{q}{\rm d}v_g(x)\right)^{1/q}\leq
			C
			\left(\int_{M} d_g^{-2a}(x_0,x)|\nabla_g f(x)|^2{\rm d}v_g(x)\right)^{1/2}},
	\end{equation}
then for every $x\in M$ and $r>0,$ we have
\begin{equation}\label{kovetkeztetes}
	{\rm Vol}_g(B_x(r))\geq \left(\frac{K_{a,b}}{C}\right)^\frac{n}{a+1-b}\omega_nr^n.
\end{equation}
 \end{theorem}

\begin{remark}\rm 
	(a) We first notice that the assumption $C\geq K_{a,b}$ is superfluous in Theorem \ref{thm-docarmo}. Indeed, the validity of \eqref{egyenlet-ujra} together with a careful  analysis based on local charts will imply $C\geq K_{a,b}$; a similar argument can be found e.g.\ in Hebey \cite{Hebey} and Krist\'aly \cite{Kristaly-JMPA}.  
	
	(b) We mention that the original proof from \cite{doCarmo-Xia} is based on a lengthy local and global   comparison of an ordinary differential inequality to an ODE. Inspired by our proof of sharpness in \S \ref{section-3}, we provide an alternative, simplified way to prove Theorem \ref{thm-docarmo}.
\end{remark}

{\it Proof of Theorem \ref{thm-docarmo}.} Let $(M,g)$ be an  $n$-dimensional Riemannian manifold $(n\geq 3)$ as in the statement, and $x_0\in M$ be fixed.  
	Let $s,t>0$ and $r\in \mathbb R$ be some constants such that
	$$st>n+r>0.$$
	A similar argument as in Proposition \ref{AT-estimates-proposition}, based on the layer cake representation and the Bishop-Gromov volume comparison principle, shows that the functional 
	$$
		\mathcal K(\lambda,r,t,s)=\int_Md_g^{r}(x_0,x)(\lambda+d_g^{t}(x_0,x))^{-s}{\rm d}v_g(x)
$$
	is well-defined for every $\lambda>0$ and the following asymptotic property holds: 
	\begin{equation} \label{eq-AT-1-000} 
	\lim_{\lambda\to \infty} \lambda^{s-\frac{n+r}{t}} \mathcal K(\lambda,r,t,s)= \frac{n}{t}\omega_n{\sf AVR}_g\frac{\Gamma\left(\frac{n+r}{t}\right)\Gamma\left(s-\frac{n+r}{t}\right)}{\Gamma(s)}.
\end{equation} 

Let $\lambda>0$. Since the Talentian bubble $M\ni x\mapsto \left(\lambda+d_g(x_0,x)^{2-bq+2a}\right)^{-\frac{n-2a-2}{2-bq+2a}}$ can be approximated by functions belonging to $C_0^\infty(M)$, we may use it as a test function in 
\eqref{egyenlet-ujra}, obtaining
$$	\displaystyle  { \left(\mathcal K\left(\lambda,-bq,2-bq+2a,q\frac{n-2a-2}{2-bq+2a}\right)\right)^{1/q}}$$ $$\leq
	C(n-2a-2)
	\left(\mathcal K\left(\lambda,2a+2-2bq,2-bq+2a,\frac{2(n-bq)}{2-bq+2a}\right)\right)^{1/2}.$$
According to the assumptions from \eqref{a-b}, the above expressions are well defined.  Let us multiply the latter inequality by 
$\lambda^\frac{n-2+2(b-a)}{4(a+1-b)}>0$ and take  the limit $\lambda\to \infty$; then, by \eqref{eq-AT-1-000} and trivial identities between $a,b,n$ and $q$, we obtain that 
$$	\displaystyle  { \left(\frac{n}{2-bq+2a}\omega_n{\sf AVR}_g\frac{\Gamma^2\left(\frac{n}{2(a+1-b)}\right)}{\Gamma(\frac{n}{a+1-b})}\right)^{1/q}}$$ $$\leq
C(n-2a-2)
\left(\frac{n}{2-bq+2a}\omega_n{\sf AVR}_g\frac{\Gamma\left(\frac{n}{2(a+1-b)}+1\right)\Gamma\left(\frac{n}{2(a+1-b)}-1\right)}{\Gamma(\frac{n}{a+1-b})}\right)^{1/2}.$$
Note that we necessarily have ${\sf AVR}_g>0$; otherwise, since $q>2$ (which is equivalent to $b<a+1$), the latter inequality would give a contradiction. Reorganizing the latter inequality, and taking into account the form of $K_{a,b}$, one obtains
 that
 $$K_{a,b}{\sf AVR}_g^{-\frac{a+1-b}{n}}\leq C.$$
 By the definition of ${\sf AVR}_g$ and the non-increasing property of  $r\mapsto \frac{{\rm Vol}_g(B_x(r))}{r^n}$, $r>0$ (cf. Bishop-Gromov comparison principle), we have for every $x\in M$ and $r>0$ that
 $$\frac{{\rm Vol}_g(B_x(r))}{\omega_nr^n}\geq {\sf AVR}_g \geq \left(\frac{K_{a,b}}{C}\right)^\frac{n}{a+1-b},$$
 as desired in \eqref{kovetkeztetes}.  
  \hfill $\square$

\begin{remark}\rm 
 The fine asymptotic properties of the Talantian bubble, described in \eqref{eq-AT-1-000}, essentially simplifies not only the proof of do Carmo and Xia \cite{doCarmo-Xia}, but also further related results, as those from Krist\'aly \cite{Kristaly-Calculus}, Krist\'aly and Ohta \cite{Kristaly-Ohta}, Mao \cite{Mao}, Tokura, Adriano and  Xia \cite{WLC}, Xia \cite{Xia, Xia2}, and references therein.  
\end{remark}

\noindent \textbf{Acknowledgements.} The author would like to express his gratitude to Zolt\'an M. Balogh, Alessio Figalli and Andrea Mondino for helpful conversations on various topics related
to the paper. He also thanks the Referee 
for his/her thoughtful comments and efforts towards improving the presentation of the manuscript.
\\

\noindent \textbf{Data availibility}. 
No data, models, or code were generated or used during the study.

\end{document}